\documentclass[a4paper,11pt]{article}

\usepackage{amsmath}
\usepackage{amssymb}
\usepackage{graphicx}
\usepackage{float}
\usepackage{ulem}
\usepackage{enumerate}
\usepackage{comment}
\usepackage{here}
\usepackage{fancybox}
\usepackage[all]{xy}

\usepackage{tabularx}
   \newcolumntype{C}{>{\centering\arraybackslash}X}
   \newcolumntype{L}{>{\raggedright\arraybackslash}X}
   \newcolumntype{R}{>{\raggedleft\arraybackslash}X}
   \newcolumntype{D}{>{\scriptsize\centering\arraybackslash}X}

\newtheorem{theorem}{Theorem}
\newtheorem{reftheorem}{Theorem}

\newtheorem{claim}{Claim}[section]
\newtheorem{subclaim}{Subclaim}[claim]
\newtheorem{lemma}{Lemma}

\newtheorem{refconjecture}[reftheorem]{Conjecture}

\newtheorem{problem}{Problem}
\newtheorem{fact}[claim]{Fact}

\renewcommand{\labelenumi}{(\theenumi)}

\newcommand{\ora}{\overrightarrow}
\newcommand{\ola}{\overleftarrow}

\newcommand{\proof}{\medbreak\noindent\textit{Proof.}\quad}
\newcommand{\qed}{{$\quad\square$\vspace{3.6mm}}}

\newcommand{\sm}{\small}
\newcommand{\scr}{\scriptsize}

\numberwithin{equation}{section}

\begin{document}

\title{Minimum degree conditions for the existence of cycles of all lengths modulo $k$ in graphs}

\author{Shuya Chiba\thanks{Applied Mathematics, Faculty of Advanced Science and Technology, 
Kumamoto University, 
2-39-1 Kurokami, Kumamoto 860-8555, Japan; 
E-mail address: \texttt{schiba@kumamoto-u.ac.jp}; 
This work was supported by JSPS KAKENHI Grant Number 17K05347}
\and 
Tomoki Yamashita\thanks{Department of Mathematics, Kindai University, 
3-4-1 Kowakae, Higashi-Osaka, Osaka 577-8502, Japan; 
E-mail address: \texttt{yamashita@math.kindai.ac.jp}; 
This work was supported by JSPS KAKENHI Grant Number 16K05262}
}

\date{}

\maketitle

\vspace{-24pt}
\begin{abstract}
Thomassen, in 1983, conjectured that 
for a positive integer $k$, 
every $2$-connected non-bipartite graph of minimum degree at least $k + 1$ 
contains cycles of all lengths modulo $k$. 
In this paper, we settle this conjecture affirmatively. 

\end{abstract}


\section{Introduction}
\label{sec:introduction}

All graphs considered in this paper are finite undirected graphs without loops or multiple edges. 
In \cite{T1983}, 
Thomassen conjectured the following.

\begin{refconjecture}[Thomassen \cite{T1983}]
\label{T1983 all length k}
For a positive integer $k$, 
every $2$-connected non-bipartite graph 
of minimum degree at least $k + 1$ 
contains cycles of all lengths modulo $k$. 
\end{refconjecture}

In 2018, 
Liu and Ma proved that this conjecture is true for all even integers $k$, 
see \cite[Theorem~1.9]{LM2018} 
(for the history and other related results to the conjecture, 
we also refer the reader to \cite{LM2018}). 
In this paper, 
we settle Conjecture~\ref{T1983 all length k}  
by showing that 
it is also true for all odd integers $k$. 
For this purpose, 
we give the following result, 
which is our main theorem. 
Here $E(G)$ denotes the edge set of a graph $G$.

\begin{theorem}
\label{Conjecture 6.2 of LM2018 for 2-connected case}
For a positive integer $k$, 
every $2$-connected graph of minimum degree at least $k + 1$ 
contains $k$ cycles $C_{1}, \dots, C_{k}$ 
such that either 
{\rm (i)} $|E(C_{i+1})| - |E(C_{i})| = 1$ for $1 \le i \le k-1$, 
or 
{\rm (ii)} $|E(C_{i+1})| - |E(C_{i})| = 2$ for $1 \le i \le k-1$. 
\end{theorem}

We here prove Conjecture~\ref{T1983 all length k} assuming Theorem~\ref{Conjecture 6.2 of LM2018 for 2-connected case} 
for the case where $k$ is odd. 
It follows from the proof that 
the condition ``non-bipartite'' in Conjecture~\ref{T1983 all length k} 
is not necessary if $k$ is odd.

\medskip
\noindent
\textit{Proof of Conjecture~\ref{T1983 all length k} for the case where $k$ is odd.} 
Let $k$ be a positive odd integer, 
and let $G$ be a $2$-connected graph 
of minimum degree at least $k + 1$ (the graph $G$ may be bipartite). 
We will show that $G$ contains cycles of all lengths modulo $k$. 
By Theorem~\ref{Conjecture 6.2 of LM2018 for 2-connected case}, 
$G$ contains $k$ cycles satisfying (i) or (ii) in Theorem~\ref{Conjecture 6.2 of LM2018 for 2-connected case}. 
If the $k$ cycles satisfy (i), then the $k$ cycles clearly have all lengths modulo $k$.  
So, suppose that 
the $k$ cycles satisfy (ii). 
Since $k$ is odd, 
we may assume that 
\begin{align*}
(l, \ l + 2, \ \dots, \ l + k-3, \ l + k-1, \ l + k+1, \ \dots, \ l + 2k-4, \ l + 2k-2)
\end{align*} 
is a sequence of the lengths of the $k$ cycles for some integer $l \ge 3$. 
Then 
it follows that 
\begin{align*}
l + 2i + 1 \equiv l + k + 2i + 1 \!\! \pmod{k} \ \textup{ for } \ 0 \le i \le \frac{k-3}{2}. 
\end{align*}
Thus the $k$ cycles have all lengths modulo $k$.
\qed

Our proof of Theorem~\ref{Conjecture 6.2 of LM2018 for 2-connected case} 
is based on the technique of Liu and Ma \cite{LM2018}. 
In the next section, 
we introduce results to prove Theorem~\ref{Conjecture 6.2 of LM2018 for 2-connected case}. 
In particular, 
we give sharp degree conditions 
for the existence of paths with specified end vertices 
whose lengths differ by one or two 
(see Theorems~\ref{improvement of LM2018 Lemma 3.1} and \ref{improvement of LM2018 Lemma 3.1 No.2} for the detail), 
which are our key results.



\section{Outline of the proof of Theorem~\ref{Conjecture 6.2 of LM2018 for 2-connected case}}
\label{sec:outline}

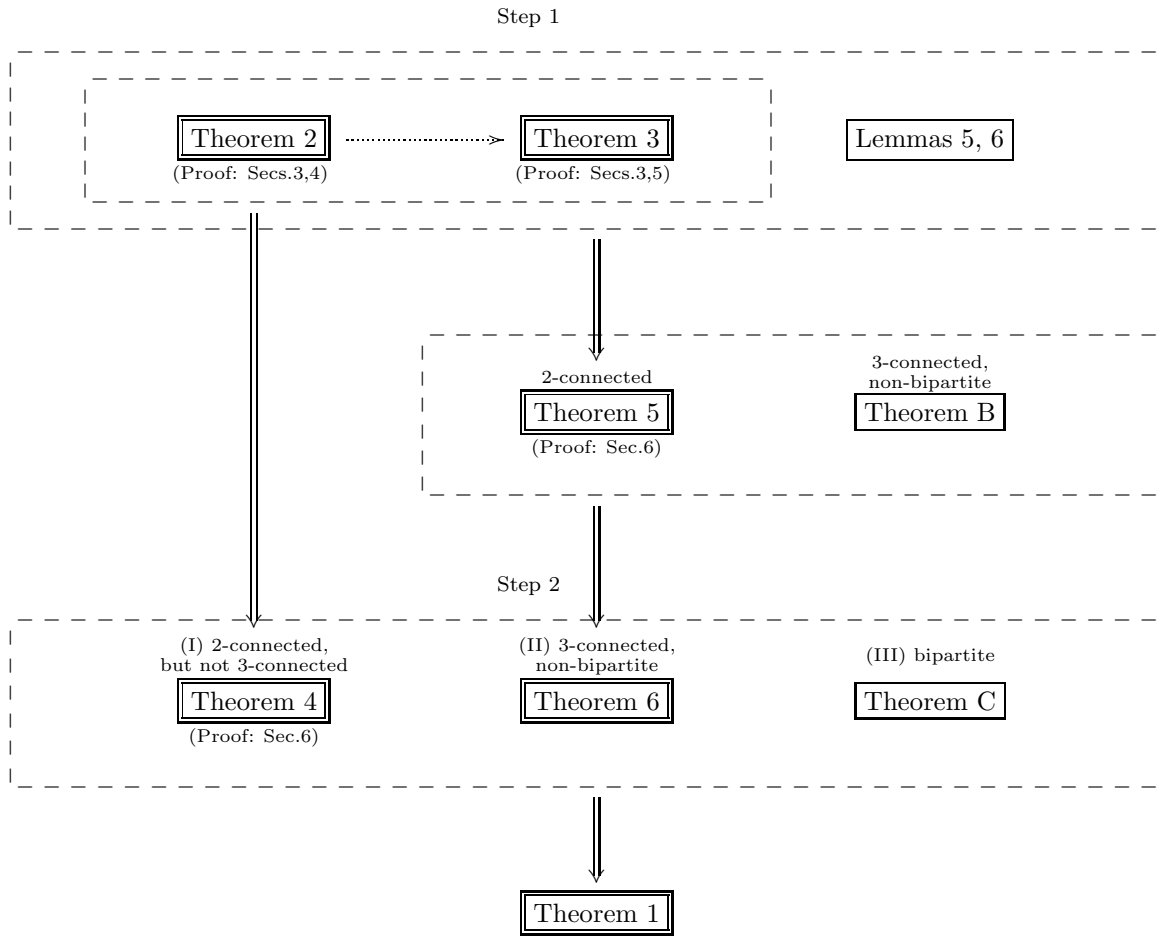
\begin{figure}[H]
{\small 
\begin{center}
\[
\xymatrix@C=22pt@R=4pt{
&%
&%
&%
&\textup{\scriptsize Step~1 \hspace{+48pt}}
&%
&%
&%
&%
\\
  \ar@{--}[rrrrrrrr]
  \ar@{--}'[dddd]
&%
&%
&%
&%
&%
&%
&%
&\ar@{--}'[dddd]
\\
&\ar@{--}[rrrr]
  \ar@{--}'[dd]
&%
&%
&%
&\ar@{--}'[dd]
&%
&%
&%
\\
&%
&\underset{ \textup{(Proof: Secs.\ref{sec:preliminaries},\ref{sec:proof of improvement of LM2018 Lemma 3.1}) } }{ \textup{ \doublebox{Theorem~\ref{improvement of LM2018 Lemma 3.1}} } }
   \ar@{.>}[rr]
&%
&\underset{ \textup{(Proof: Secs.\ref{sec:preliminaries},\ref{sec:proof of improvement of LM2018 Lemma 3.1 No.2}) } }{ \textup{ \doublebox{Theorem~\ref{improvement of LM2018 Lemma 3.1 No.2}} } }
&%
&\textup{\fbox{Lemmas~\ref{LM2018 Lemma 5.1}, \ref{lemma:consecutive cycles}}}  
&%
&%
\\
&\ar@{--}[rrrr]
&\ar@{=>}[dddddddddddd]
&%
&%
&%
&%
&%
&%
\\
  \ar@{--}[rrrrrrrr]
&%
&%
&%
&\ar@{=>}[ddddd]
&%
&%
&%
&%
\\
&%
&%
&%
&%
&%
&%
&%
&%
\\
&%
&%
&%
&%
&%
&%
&%
&%
\\
&%
&%
&%
&%
&%
&%
&%
&%
\\
&%
&%
&\ar@{--}[rrrrr]
  \ar@{--}'[dd]
&%
&%
&%
&%
&\ar@{--}'[dd]
\\
&%
&%
&%
&\overset{\substack{ \textup{$2$-connected} \\[-1.5mm] \tiny \ } }{\underset{\textup{(Proof: Sec.\ref{sec:proofs of improvement of LM2018 Lemma 4.1 and Theorem 5.2})}}{\textup{ \doublebox{Theorem~\ref{improvement of LM2018 Theorem 5.2}} }} }
&%
&\overset{ \substack{ \textup{$3$-connected,} \\ \textup{non-bipartite}} }{\textup{ \fbox{Theorem~\ref{BV1998}} }}
&%
&%
\\
&%
&%
&\ar@{--}[rrrrr]
&\ar@{=>}[ddddd]
&%
&%
&%
&%
\\
&%
&%
&%
&%
&%
&%
&%
&%
\\
&%
&%
&%
&%
&%
&%
&%
&%
\\
&%
&%
&%
&\textup{\scriptsize Step~2 \hspace{+48pt}}%
&%
&%
&%
&%
\\
  \ar@{--}[rrrrrrrr]
  \ar@{--}'[dd]
&%
&%
&%
&%
&%
&%
&%
&\ar@{--}'[dd]
\\
&%
&\overset{\substack{ \textup{(I) $2$-connected,} \\ \textup{but not $3$-connected} \\[-1.5mm] \tiny \ } }{ \underset{ \textup{(Proof: Sec.\ref{sec:proofs of improvement of LM2018 Lemma 4.1 and Theorem 5.2})} }{ \textup{ \doublebox{Theorem~\ref{improvement of LM2018 Lemma 4.1}} } } }
&%
&\overset{ \substack{ \textup{(II) $3$-connected,} \\ \textup{non-bipartite} \\[-1.5mm] \tiny \ } }{ \textup{ \doublebox{Theorem~\ref{improvement of LM2018 Theorem 1.4}} } }
&%
&\overset{ \substack{ \textup{(III) bipartite} \\[0.5mm] \tiny \ } }{ \textup{ \fbox{Theorem~\ref{LM2018 Theorem 1.2}} } }
&%
&%
\\
  \ar@{--}[rrrrrrrr]
&%
&%
&%
&\ar@{=>}[dddd]
&%
&%
&%
&%
\\
&%
&%
&%
&%
&%
&%
&%
&%
\\
&%
&%
&%
&%
&%
&%
&%
&%
\\
&%
&%
&%
&%
&%
&%
&%
&%
\\
&%
&%
&%
&\textup{\doublebox{Theorem~\ref{Conjecture 6.2 of LM2018 for 2-connected case}}}
&%
&%
&%
&%
}
\]
\end{center}
\vspace{-12pt}
\caption{The overall flowchart of the proof of Theorem~\ref{Conjecture 6.2 of LM2018 for 2-connected case}}
\label{overall flow of the proof}
}
\end{figure}

In this section, 
we introduce results for the proof of Theorem~\ref{Conjecture 6.2 of LM2018 for 2-connected case} 
according to the flowchart of Figure~\ref{overall flow of the proof}. 

\vspace{-12pt}

We require some terminology and notation. 
Let $G$ be a graph. 
The vertex set of $G$ is denoted by $V(G)$. 
For a vertex $v$ of $G$, 
$\deg_{G}(v)$ denotes the degree of $v$ in $G$, 
and let $\delta(G)$ denote the minimum degree of $G$. 
For $S \subseteq V(G)$, 
$G[S]$ denotes the subgraph of $G$ induced by $S$, 
and let $G - S = G[V(G) \setminus S]$. 
For distinct vertices $x$ and $y$ of $G$, 
$(G, x, y)$ is called a \textit{rooted graph}. 
A rooted graph $(G, x, y)$ is \textit{$2$-connected} 
if 
\begin{enumerate}
\renewcommand{\labelenumi}{\upshape{(R\arabic{enumi}})}
\item 
$G$ is a connected graph of order at least $3$ with at most two end blocks, 
and 
\item 
every end block of $G$ contains at least one of $x$ and $y$ as a non-cut vertex. 
\end{enumerate}
Note that $(G, x, y)$ is $2$-connected if and only if $G + xy$ 
(i.e., the graph obtained from $G$ by adding the edge $xy$ if $xy \notin E(G)$) is $2$-connected. 
For convenience, 
we say that 
a sequence of paths or cycles $H_{1}, H_{2}, \dots, H_{k}$ 
\begin{itemize}
\item 
\textit{have consecutive lengths} if $|E(H_{1})| \ge 2$ and $|E(H_{i+1})| - |E(H_{i})| = 1$ for $1 \le i \le k - 1$; 
\item
\textit{satisfy the length condition} if $|E(H_{1})| \ge 2$ and $|E(H_{i+1})| - |E(H_{i})| = 2$ for $1 \le i \le k - 1$; 
\item
\textit{satisfy the semi-length condition} if $|E(H_{1})| \ge 2$ and there exists an index $j$ with $1 \le j \le k - 1$, 
which is called a \textit{switch}, 
such that 
\vspace{-4pt}
\begin{align*}
&|E(H_{i+1})| - |E(H_{i})| = 2 \ \textup{for} \ 1 \le i \le j-1,\\ 
&|E(H_{j+1})| - |E(H_{j})| = 1 \ \textup{and}\\ 
&|E(H_{i+1})| - |E(H_{i})| = 2 \ \textup{for} \ j+1 \le i \le k-1. 
\end{align*}
\end{itemize}

The first step of the proof of Theorem~\ref{Conjecture 6.2 of LM2018 for 2-connected case} 
is to show the following two results concerning degree conditions 
for the existence of paths satisfying the length condition or the semi-length condition.

\begin{theorem}
\label{improvement of LM2018 Lemma 3.1}
Let $k$ be a positive integer, 
and let 
$(G, x, y)$ be a $2$-connected rooted graph. 
Suppose that 
$\deg_{G}(v) \ge 2k$ for any $v \in V(G) \setminus \{x, y\}$.
Then 
$G$ contains $k$ paths from $x$ to $y$ 
satisfying the length condition.  
\end{theorem}

\begin{theorem}
\label{improvement of LM2018 Lemma 3.1 No.2}
Let $k$ be a positive integer, 
and let 
$(G, x, y)$ be a $2$-connected rooted graph. 
Suppose that 
$\deg_{G}(v) \ge 2k-1$ for any $v \in V(G) \setminus \{x, y\}$.
Then 
$G$ contains $k$ paths from $x$ to $y$ 
satisfying the length condition or the semi-length condition.  
\end{theorem}

The complete graphs of orders $2k$ and $2k-1$ show the sharpness 
of the degree conditions in Theorems~\ref{improvement of LM2018 Lemma 3.1} and \ref{improvement of LM2018 Lemma 3.1 No.2}, respectively. 
Theorem~\ref{improvement of LM2018 Lemma 3.1} is an improvement of \cite[Lemma~3.1]{LM2018}.

In Section~\ref{sec:preliminaries}, 
we prepare lemmas for 
the proofs of Theorems~\ref{improvement of LM2018 Lemma 3.1} and \ref{improvement of LM2018 Lemma 3.1 No.2}. 
We will use Theorem~\ref{improvement of LM2018 Lemma 3.1} 
in a part of the proof of Theorem~\ref{improvement of LM2018 Lemma 3.1 No.2} 
(see also Figure~\ref{overall flow of the proof}). 
So, we first prove Theorem~\ref{improvement of LM2018 Lemma 3.1} in Section~\ref{sec:proof of improvement of LM2018 Lemma 3.1} 
and then we give the proof of Theorem~\ref{improvement of LM2018 Lemma 3.1 No.2} in Section~\ref{sec:proof of improvement of LM2018 Lemma 3.1 No.2}.

\medskip

In the second step, 
we divide the proof of Theorem~\ref{Conjecture 6.2 of LM2018 for 2-connected case} 
into three cases 
according as a graph 
is 
(I) $2$-connected, but not $3$-connected, 
(II) $3$-connected and non-bipartite, 
or 
(III) bipartite.

For the case (I), 
we show the following theorem 
by using Theorems~\ref{improvement of LM2018 Lemma 3.1} and \ref{improvement of LM2018 Lemma 3.1 No.2}, 
which is an improvement of \cite[Lemma~4.1]{LM2018}. 
We give the proof in Section~\ref{sec:proofs of improvement of LM2018 Lemma 4.1 and Theorem 5.2}.

\begin{theorem}
\label{improvement of LM2018 Lemma 4.1}
Let $k$ be a positive integer, 
and let 
$G$ be a $2$-connected but not $3$-connected graph. 
If $\delta(G) \ge k + 1$, 
then 
$G$ contains $k$ cycles satisfying the length condition. 
\end{theorem}

To show the case (II), 
in Section~\ref{sec:proofs of improvement of LM2018 Lemma 4.1 and Theorem 5.2}, 
we also prove the following theorem 
by using Theorems~\ref{improvement of LM2018 Lemma 3.1}, \ref{improvement of LM2018 Lemma 3.1 No.2} 
and additional lemmas. 
Here, 
for a cycle $C$ in a connected graph $G$, 
$C$ is said to be \textit{non-separating} if $G-V(C)$ is connected.

\begin{theorem}
\label{improvement of LM2018 Theorem 5.2}
Let $k$ be a positive integer, 
and let 
$G$ be a $2$-connected graph containing a non-separating induced odd cycle. 
If $\delta(G) \ge k + 1$, 
then 
$G$ contains $k$ cycles, 
which have consecutive lengths or satisfy the length condition.
\end{theorem}

The following result is known for the existence of 
a non-separating induced odd cycle.

\begin{reftheorem}[see the proof of Theorem~2 in \cite{BV1998}]
\label{BV1998}
Every $3$-connected non-bipartite graph contains a non-separating induced odd cycle. 
\end{reftheorem}

Combining Theorems~\ref{improvement of LM2018 Theorem 5.2} and \ref{BV1998},
we can obtain the following theorem for the case (II).

\begin{theorem}
\label{improvement of LM2018 Theorem 1.4}
Let $k$ be a positive integer, 
and let 
$G$ be a $3$-connected non-bipartite graph. 
If $\delta(G) \ge k + 1$, 
then 
$G$ contains $k$ cycles, 
which have consecutive lengths or satisfy the length condition.
\end{theorem}

On the other hand, 
for the case (III), 
the following theorem is proved by Liu and Ma.

\begin{reftheorem}[\textup{\cite[Theorem~1.2]{LM2018}}]
\label{LM2018 Theorem 1.2}
Let $k$ be a positive integer, 
and let 
$G$ be a bipartite graph. 
If $\delta(G) \ge k + 1$, then 
$G$ contains $k$ cycles satisfying the length condition.
\end{reftheorem}

Consequently, 
we can obtain Theorem~\ref{Conjecture 6.2 of LM2018 for 2-connected case} 
by Theorems~\ref{improvement of LM2018 Lemma 4.1}, \ref{improvement of LM2018 Theorem 1.4} and \ref{LM2018 Theorem 1.2}. 
In Section~\ref{sec:concluding remarks}, 
we give some remarks on 
the minimum degree and the connectivity conditions in Theorem~\ref{Conjecture 6.2 of LM2018 for 2-connected case}.

\medskip

In the rest of this section, we prepare terminology and notation 
which will be used in the subsequent sections. 
Let $G$ be a graph. 
We denote by $N_{G}(v)$ the neighborhood of a vertex $v$ in $G$. 
For $S \subseteq V(G)$, 
we define $N_{G}(S) = \left( \bigcup_{v \in S}N_{G}(v) \right) \setminus S$. 
For $S, T \subseteq V(G)$ with $S \cap T = \emptyset$, 
$E_{G}(S, T)$ denotes the set of edges of $G$ between $S$ and $T$, 
and let $e_{G}(S, T) = |E_{G}(S, T)|$. 
Furthermore, 
$G[S, T]$ is the graph defined by 
$V(G[S, T])=S \cup T$ and 
$E(G[S, T])=E_{G}(S, T)$.
Note that $G[S, T]$ is a bipartite subgraph of $G$ 
with partite sets $S$ and $T$, 
and we always assume that 
$G[S, T]$ is such a bipartite graph. 
For a rooted graph $(G, x, y)$, 
we define $\delta(G, x, y) = \min \{\deg_{G}(v) : v \in V(G) \setminus \{x, y\}\}$. 
In the rest of this paper, 
we often denote the singleton set $\{v\}$ by $v$, 
and 
we often identify a subgraph $H$ of $G$ with its vertex set $V(H)$.

For $S \subseteq V(G)$, 
a path in $G$ is an $S$-\textit{path} 
if it begins and ends in $S$, 
and 
none of its internal vertices are contained in $S$. 
For $S, T \subseteq V(G)$ with $S \cap T = \emptyset$, 
a path in $G$ is an $(S, T)$-path 
if 
one end vertex of the path belongs to $S$, another end vertex belongs to $T$, 
and the internal vertices do not belong to $S \cup T$. 
We write a path or a cycle $P$ with a given orientation as $\ora{P}$. 
If there exists no fear of confusion, 
we abbreviate $\ora{P}$ by $P$. 
Let $\ora{P}$ be an oriented path (or cycle). 
For $u \in V(P)$, 
the $h$-th successor and the $h$-th predecessor of $u$ 
on $\ora{P}$ (if exist) 
is denoted by $u^{+h}$ and $u^{-h}$, respectively, 
and we let $u^{+} = u^{+1}$ and $u^{-} = u^{-1}$. 
For $u, v \in V(P)$, 
the path from $u$ to $v$ along $\ora{P}$ (if exist) 
is denoted by $u\ora{P}v$. 
The reverse sequence of $u \ora{P} v$ is denoted by $v \ola{P} u$. 
In the rest of this paper, 
if $\ora{P}$ is an $(S, T)$-path in $G$, 
we always assume that 
the orientation of $P$ is given
from the end vertex belonging to $S$
to the end vertex belonging to $T$
along the edges of $P$.



\section{Preliminaries for the proofs of Theorems~\ref{improvement of LM2018 Lemma 3.1} and \ref{improvement of LM2018 Lemma 3.1 No.2}}
\label{sec:preliminaries}

In this section, 
we introduce the concept of a core which was used in the argument of \cite{LM2018} 
and give some lemmas 
for the proofs of Theorems~\ref{improvement of LM2018 Lemma 3.1} and \ref{improvement of LM2018 Lemma 3.1 No.2}. 

Let $x$ and $y$ be two distinct vertices of a graph $G$. 
For a bipartite subgraph $H = G[S, T]$ of $G$ and an integer $l$, 
$H$ is called an $l$-\textit{core with respect to $(x, y)$} 
if 
\begin{enumerate}
\renewcommand{\labelenumi}{\upshape{(C\arabic{enumi}})}
\item 
$H$ is complete bipartite and $|T| \ge |S| = l + 1 \ge 2$, 
\item 
$x \in S$ and $y \notin V(H)$, 
\item 
$e_{G}(v, S) \le l$ for $v \in V(G) \setminus (V(H) \cup \{y\})$, and 
\item 
$e_{G}(v, T \setminus \{v\}) \le l + 1$ for $v \in V(G) \setminus (S \cup \{y\})$. 
\end{enumerate}
In the rest of this section, 
we fix the following notation. 
Let $(G, x, y)$ be a $2$-connected rooted graph, 
and let $H = G[S, T]$ be an $l$-core of $G$ with respect to $(x, y)$. 
Furthermore, 
let $C$ be the component of $G- V(H)$ such that $y \in V(C)$. 
Since $E_{G}(H-x, C) \neq \emptyset$ by (R2), 
the following two lemmas (Lemmas~\ref{lemma:l >= k or l >= k-1} and \ref{lemma:l >= k-1, E(T) is not empty}) easily follows from (C1). 
So, we omit the proof.

\begin{lemma}
\label{lemma:l >= k or l >= k-1} 
If either 
{\rm (i)}~$l \ge k$ 
or 
{\rm (ii)}~$l = k -1$ and $E_{G}(T, C) \neq \emptyset$, 
then $G$ contains $k$ $(x, y)$-paths satisfying the length condition.
\end{lemma}

\begin{lemma}
\label{lemma:l >= k-1, E(T) is not empty} 
If $l = k-1$, 
$E_{G}(S \setminus \{x\}, C) \neq \emptyset$ 
and 
$E(G[T]) \neq \emptyset$, 
then $G$ contains $k$ $(x, y)$-paths satisfying the semi-length condition. 
\end{lemma}

The following two lemmas 
(Lemmas~\ref{lemma:paths satisfying the length condition (1)} and \ref{lemma:paths satisfying the length condition (2)}) 
are proved by Liu and Ma, see \cite[Lemmas~2.7 and 2.11]{LM2018}. 
Note that the argument in \cite{LM2018} can work for paths satisfying the length condition 
but also for paths satisfying the semi-length condition.

\begin{lemma}
\label{lemma:paths satisfying the length condition (1)} 
Let $s$ be a vertex of $S \setminus \{x\}$ such that 
$E_{G}(s, C) \neq \emptyset$. 
If one of the following {\rm (i)--(iii)} holds, 
then $G$ contains $k$ $(x, y)$-paths satisfying the length condition 
(resp., the semi-length condition). 
\begin{enumerate}[{\upshape(i)}]
\item 
\label{T to T}
$G-V(C)$ contains $k-l+1$ $T$-paths internally disjoint from $V(H)$ 
and satisfying the length condition 
(resp., the semi-length condition) \textup{{\bf \cite[Lemma~2.7-2]{LM2018}}}. 

\item 
\label{T to {x,s}}
$G-V(C)$ contains $k-l+1$ $(T, \{x, s\})$-paths internally disjoint from $V(H)$ 
and satisfying the length condition 
(resp., the semi-length condition) \textup{{\bf \cite[Lemma~2.7-4]{LM2018}}}. 

\item 
\label{T to S-{x,s}}
$G-V(C)$ contains $k - l + 2$ $(T, S \setminus \{x, s\})$-paths
internally disjoint from $V(H)$ 
and satisfying the length condition 
(resp., the semi-length condition) \textup{{\bf \cite[Lemma~2.7-3]{LM2018}}}. 

\end{enumerate}
\end{lemma}
\begin{lemma}[\textup{\cite[Lemma~2.11]{LM2018}}]
\label{lemma:paths satisfying the length condition (2)} 
If one of the following {\rm (i)} and {\rm (ii)} holds, 
then $G$ contains $k$ $(x, y)$-paths satisfying the length condition 
(resp., the semi-length condition). 
\begin{enumerate}[{\upshape(i)}]
\item
\label{T to y}
$G$ contains $k-l$ $(T, y)$-paths 
internally disjoint from $V(H)$ 
and satisfying 
the length condition 
(resp., the semi-length condition). 

\item
\label{S - x to y}
$G$ contains $k - l + 1$ $(S \setminus \{x\}, y)$-paths
internally disjoint from $V(H)$ 
and satisfying 
the length condition 
(resp., the semi-length condition). 

\end{enumerate}
\end{lemma}



\section{Proof of Theorem~\ref{improvement of LM2018 Lemma 3.1}}
\label{sec:proof of improvement of LM2018 Lemma 3.1}

\noindent
\textit{Proof of Theorem~\ref{improvement of LM2018 Lemma 3.1}.}~We prove it by induction on $|V(G)| + |E(G)|$.
Let $(G, x, y)$ be a minimum counterexample with respect to $|V(G)| + |E(G)|$. 
If $k = 1$, then by (R1) and (R2), 
we can easily see that $G$ contains an $(x, y)$-path of length at least $2$, 
a contradiction. 
Thus $k \ge 2$. 
Since $\delta(G, x, y) \ge 2k$, 
this implies that $|G| \ge 5$. 
By symmetry, we may assume that 
$\deg_{G}(x) \le \deg_{G}(y)$.

\begin{claim}
\label{claim:2-connected} 
$G$ is $2$-connected. 
\end{claim}
\proof 
Suppose that $G$ is not $2$-connected. 
Then by (R1), 
$G$ has a cut vertex $c$ 
and $G-c$ has exactly two components $C_{1}$ and $C_{2}$. 
Since $|G| \ge 5 \ (\ge 4)$, we may assume that $|C_{1}| \ge 2$.
Let $G_{i} = G[V(C_{i}) \cup \{c\}]$ for $i \in \{1, 2\}$. 
Then by (R2), 
for some two distinct vertices $x', y' \in \{x, y\}$, 
(i)~$(G_{1}, x', c)$ is a 2-connected rooted graph 
such that 
$\delta(G_{1}, x', c) \ge \delta(G, x, y) \ge 2k$, 
and (ii)~$y'$ is contained in a block of $G_{2}$. 
Hence by the induction hypothesis, 
$G_{1}$ contains $k$ $(x', c)$-paths $\ora{P_{1}}, \dots, \ora{P_{k}}$ satisfying the length condition. 
Then $x'\ora{P_{1}}c \ora{P} y', \dots, x'\ora{P_{k}}c \ora{P} y'$ 
are $k$ $(x, y)$-paths in $G$ satisfying the length condition, 
where 
$\ora{P}$ denotes the $(c, y')$-path in $G_{2}$, 
a contradiction. 
\qed

\begin{claim}
\label{claim:xy is not an edge} 
$xy \notin E(G)$. 
\end{claim}
\proof 
If $xy \in E(G)$, 
then 
by Claim~\ref{claim:2-connected}, 
$(G - xy, x, y)$ is a $2$-connected rooted graph 
such that 
$\delta(G-xy, x, y) = \delta(G, x, y) \ge 2k$, 
and hence 
the induction hypothesis yields that 
$G -xy$ (and also $G$) contains $k$ $(x, y)$-paths 
satisfying the length condition, a contradiction. 
\qed

\begin{enumerate}[{\textup{{\bf Case~\arabic{enumi}.}}}]
\setcounter{enumi}{0}
\item 
$G-y$ does not contain a cycle of length $4$ passing through $x$. 
\end{enumerate}

Since $xy \notin E(G)$ by Claim~\ref{claim:xy is not an edge}, 
in this case, 
we have 
\begin{align}
\label{e(v, N_{G}(x))}
e_{G}(v, N_{G}(x) \setminus \{v\}) \le 1 \textup{ for } v \in V(G) \setminus \{x, y\}. 
\end{align}
Let $G^{*}$ be the graph obtained from $G$ 
by contracting the subgraph induced by $N_{G}(x) \cup \{x\}$ 
into a single vertex $x^{*}$ and then removing multiple edges. 
Then by (\ref{e(v, N_{G}(x))}), 
\begin{align}
\label{minimum degree of G^{*}}
\deg_{G^{*}}(v) = \deg_{G}(v) \textup{ for } v \in V(G^{*}) \setminus \{x^{*}, y\}. 
\end{align}
By (\ref{e(v, N_{G}(x))}) and since $\delta(G, x, y) \ge 2k \ge 4$, 
we have $|G^{*}| \ge 3$. 
Hence 
by Claim~\ref{claim:2-connected},  
if $G^{*}$ is not $2$-connected, 
then $x^{*}$ is the unique cut vertex of $G^{*}$ 
and each block of $G^{*}$ is an end block containing $x^{*}$. 
Now, 
let $B^{*}$ be the block of $G^{*}$ which contains $y$ 
if $G^{*}$ is not $2$-connected; 
otherwise, let $B^{*} = G^{*}$.

Assume that 
$(B^{*}, x^{*}, y)$ is $2$-connected. 
Since $\delta(B^{*}, x^{*}, y) \ge \delta(G, x, y) \ge 2k$ 
by (\ref{minimum degree of G^{*}}), 
it follows from the induction hypothesis that 
$B^{*}$ contains $k$ $(x^{*}, y)$-paths $\ora{P_{1}}, \dots, \ora{P_{k}}$ 
satisfying the length condition. 
By the definition of $G^{*}$ 
and since $xy \notin E(G)$, 
we see that 
for each $i$ with $1 \le i \le k$, 
there is a vertex $u_{i}$ of $N_{G}(x)$ 
such that $u_{i} u'_{i} \in E(G)$, 
where $u'_{i}$ is the successor of $x^{*}$ along $\ora{P_{i}}$. 
Therefore, 
$x u_{1} u'_{1} \ora{P_{1}} y, \dots, x u_{k} u'_{k} \ora{P_{k}} y$ 
are $k$ $(x, y)$-paths in $G$ satisfying the length condition, a contradiction. 
Thus 
$(B^{*}, x^{*}, y)$ is not $2$-connected.

Then we have $V(B^{*}) = \{x^{*}, y\}$. 
This together with the definition of $G^{*}$ and Claim~\ref{claim:xy is not an edge} 
implies that $N_{G}(y) \subseteq N_{G}(x)$. 
Since $\deg_{G}(x) \le \deg_{G}(y)$, 
this yields that $N_{G}(x) = N_{G}(y)$. 
Since $|G^{*}| \ge 3$, the definition of $G^{*}$ also implies that 
there is a component $C$ of $G - (N_{G}(x) \cup \{x, y\})$. 
Then $N_{G}(C) \subseteq N_{G}(x)$ 
and 
$G[V(C) \cup N_{G}(C) \cup \{x\}]$ is $2$-connected, 
since $G$ is $2$-connected and $N_{G}(x) = N_{G}(y)$. 
Hence $N_{G}(C) \ (\subseteq N_{G}(x))$ can be partitioned 
into two vertex-disjoint non-empty sets $S$ and $T$ 
so that 
the graph $C^{*}$ defined as follows is $2$-connected: 
\begin{align*}
V(C^{*}) 
&= V(C) \cup \{x, S, T\} \textup{ and} \\
E(C^{*})
&= E(C) \cup \{xS, xT\} 
\cup \{vS : v \in V(C), E_{G}(v, S) \neq \emptyset\} \cup \{vT : v \in V(C), E_{G}(v, T) \neq \emptyset\}. 
\end{align*}
Then 
by the definition of $C^{*}$, 
it follows that
$(C^{*} - x, S, T)$ is a $2$-connected rooted graph 
and $\deg_{C^{*} - x}(v) = \deg_{G}(v)$ for $v \in V(C^{*}) \setminus \{x, S, T\}$; 
thus $\delta(C^{*} - x, S, T) \ge \delta(G, x, y) \ge 2k$. 
By the induction hypothesis, 
$C^{*}-x$ contains $k$ $(S, T)$-paths $\ora{P_{1}}, \dots, \ora{P_{k}}$
satisfying the length condition.
Note that each $P_{i}$ has order at least $3$, 
and each $P_{i} - \{S, T\}$ is contained in $C$. 
For each $i$ with $1 \le i \le k$, 
let $s_{i}$ and $t_{i}$ be vertices in $S$ and $T$, respectively, 
such that 
$s_{i} s'_{i}, t_{i} t'_{i} \in E(G)$, 
where $s'_{i}$ and $t'_{i}$ are the successor of $S$
and the predecessor of $T$ along $\ora{P_{i}}$, respectively. 
Then 
$x s_{1} s'_{1} \ora{P_{1}} t'_{1} t_{1} y , \dots, x s_{k} s'_{k} \ora{P_{k}} t'_{k} t_{k} y$ 
are $k$ $(x, y)$-paths in $G$ satisfying the length condition, a contradiction.

This completes the proof of Case~1.

\begin{enumerate}[{\textup{{\bf Case~\arabic{enumi}.}}}]
\setcounter{enumi}{1}
\item 
$G-y$ contains a cycle of length $4$ passing through $x$. 
\end{enumerate}

By the assumption of Case~2, 
$G$ contains a bipartite subgraph $H = G[S, T]$ 
such that 
$H$ is complete bipartite, $|T| \ge |S| =: l + 1 \ge 2$, 
$x \in S$, $y \notin V(H)$ (i.e., $H$ satisfies (C1) and (C2)). 
Let $C$ be the component of $G-V(H)$ 
such that $y \in V(C)$. 
Choose $H$ so that 
\begin{enumerate}[{\upshape(a)}]
\item $|S|$ is maximum, 
\item $|T|$ is maximal, subject to (a), 
\item $|C|$ is maximum, subject to (a) and (b), and 
\item $|N_{G}(C) \cap S|$ is minimum, subject to (a), (b) and (c). 
\end{enumerate}
Then by the choices (a) and (b), we can obtain the following.

\begin{claim}
\label{claim:l-core} 
$H$ is an $l$-core of $G$ with respect to $(x, y)$. 
\end{claim}
\proof 
Since $H$ satisfies (C1) and (C2), 
it suffices to show that $H$ also satisfies (C3) and (C4).

We first show (C4). 
Suppose that there exists a vertex $v \in V(G) \setminus (S \cup \{y\})$ 
such that 
$e_{G}(v, T \setminus \{v\}) \ge l + 2$. 
Let $S' = S \cup \{v\}$ 
and 
$T' = \{v' \in V(G) \setminus (S' \cup \{y\}) : S' \subseteq N_{G}(v')\}$. 
Note that 
$N_{G}(v) \cap T \subseteq T'$. 
Hence 
$G[S', T']$ is a complete bipartite subgraph of $G$ 
such that 
$|T'| \ge e_{G}(v, T \setminus \{v\}) \ge l + 2 = |S'| > |S|$, 
$x \in S'$ and $y \notin V(G[S', T'])$, 
which contradicts the choice (a).

We next show (C3). 
Suppose that there exists a vertex $v \in V(G) \setminus (V(H) \cup \{y\})$ 
such that 
$e_{G}(v, S) \ge l + 1$. 
Since $|S| = l + 1$, 
this implies that 
$S \subseteq N_{G}(v)$. 
Hence, 
$G[S, T \cup \{v\}]$ 
is a complete bipartite subgraph of $G$ 
such that 
$|T \cup \{v\}| > |T| \ge |S|$, $x \in S$ and $y \notin V(G[S, T \cup \{v\}])$
which contradicts the choice (b). 
\qed

By the choices (c) and (d), we can obtain the following.

\begin{claim}
\label{claim:e_{G}(v, T) <= l} 
If $E_{G}(S \setminus \{x\}, C) \neq \emptyset$, 
then 
{\rm (i)}~$e_{G}(v, T) \le l$ for $v \in V(G) \setminus V(H \cup C)$, 
and 
{\rm (ii)}~$e_{G}(v, T \setminus \{v\}) \le l$  or $E_{G}(v, C) \neq \emptyset$ for $v \in T$. 
\end{claim}
\proof 
Assume that $E_{G}(S \setminus \{x\}, C) \neq \emptyset$, 
and let $s$ be a vertex of $S \setminus \{x\}$ such that 
$E_{G}(s, C) \neq \emptyset$. 
We now 
suppose that there is a vertex $v$ of $V(G) \setminus (S \cup V(C))$ 
such that $e_{G}(v, T \setminus \{v\}) \ge l + 1$. 
Further, if $v \in T$, then we suppose that 
$E_{G}(v, C) = \emptyset$. 
Note that 
$E_{G}(v, C) = \emptyset$ also holds 
in the case of $v \in V(G) \setminus V(H \cup C)$, 
since $C$ is a component of $G-V(H)$.

Let 
$S' = (S \setminus \{s\}) \cup \{v\}$ 
and 
$T' = \{v' \in V(G) \setminus (S' \cup \{y\}) : S' \subseteq N_{G}(v')\}$, 
and let 
$H' = G[S', T']$. 
Since 
$e_{G}(v, T \setminus \{v\}) \ge l + 1$, 
it follows from 
the definitions of $S', T'$ and $H'$ 
that 
$H'$ is a complete bipartite subgraph of $G$ 
such that 
$|T'| \ge e_{G}(v, T \setminus \{ v \}) \ge l + 1 = |S'| = |S|$, 
$x \in S'$ and $y \notin V(H')$. 
In particular, $H'$ satisfies (a) and (b).

Let 
$C'$ be the component of $G-V(H')$ such that $y \in V(C')$. 
Note that 
$V(C) \subseteq V(C') \subseteq V(G) \setminus V(H')$ 
since $S' = (S \setminus \{s\}) \cup \{v\}$ 
and 
$E_{G}(v, C) = \emptyset$. 
Hence 
the choice (c) yields that 
$C' = C$. 
In particular, $H'$ also satisfies (c). 
But then, since $s \in N_{G}(C)$ and $v \notin N_{G}(C)$, 
we have 
$$| N_{G}(C') \cap S' | 
= | N_{G}(C) \cap \big( (S \setminus \{s\}) \cup \{v\} \big) | 
= | N_{G}(C) \cap S |-1,$$ 
which contradicts the choice (d). 
Thus (i) and (ii) hold. 
\qed

Note that $l \le k - 1$
by Lemma~\ref{lemma:l >= k or l >= k-1} and Claim~\ref{claim:l-core}.

\begin{claim}
\label{claim:neighbor of T} 
If $E_{G}(S \setminus \{x\}, C) \neq \emptyset$, 
then $N_{G}(T) \subseteq V(H \cup C)$. 
\end{claim}
\proof
Assume that $E_{G}(S \setminus \{x\}, C) \neq \emptyset$, 
and let $s$ be a vertex of $S \setminus \{x\}$ such that 
$E_{G}(s, C) \neq \emptyset$. 
Then we can take an $(s, y)$-path $\ora{P}$ in $G[V(C) \cup \{s\}]$.

We now suppose that $N_{G}(T) \not \subseteq V(H \cup C)$, 
i.e., 
there exists a component $D$ of $G - V(H)$ 
such that $D \neq C$ and $E_{G}(T, D) \neq \emptyset$.

\begin{subclaim}
\label{subclaim:E_{G}(B-b, T cup x, s)} 
Let $B$ be an end block of $D$, 
and let $b$ be a cut vertex of $D$ which is contained in $B$. 
Then $E_{G}(B-b, T \cup \{x, s\}) \neq \emptyset$.  
\end{subclaim}
\proof 
Suppose that $E_{G}(B-b, T \cup \{x, s\}) = \emptyset$. 
Since $G$ is $2$-connected, we have $E_{G}(B-b, S \setminus \{x, s\}) \neq \emptyset$. 
In particular, $S \setminus \{x, s\} \neq \emptyset$, that is, $\l \ge 2$. 
We define the graph $B^{*}$ as follows: 
\begin{align*}
V(B^{*}) 
&= V(B) \cup \{S^{*}\} 
\textup{ and }
E(B^{*})
= E(B) \cup \{v S^{*} : v \in V(B), E_{G}(v, S \setminus \{x, s\}) \neq \emptyset\}. 
\end{align*}
Then 
$(B^{*}, S^{*}, b)$ is a $2$-connected rooted graph. 
Since $l \ge 2$, it also follows that 
for each $v \in V(B^{*}) \setminus \{S^{*}, b\}$, 
\begin{align*}
\deg_{B^{*}}(v) \ge 
\left \{
\begin{array}{ll} 
\deg_{G}(v) \ge 2k  & \textup{(if $E_{G}(v, S \setminus \{x, s\}) = \emptyset$)} \\[3mm]
\deg_{G}(v) - (l-1) + 1 \ge 2 (k- l+2) & \textup{(if $E_{G}(v, S \setminus \{x, s\}) \neq \emptyset$)}
\end{array}
\right.,  
\end{align*}
and thus $\delta(B^{*}, S^{*}, b) \ge 2 ( k- l+2 )$. 
By the induction hypothesis, 
$B^{*}$ contains $k- l+2$ $(S^{*}, b)$-paths satisfying the length condition. 
Then it follows from the definition of $B^{*}$ that 
$G-V(C)$ contains $k- l+2$ $(S \setminus \{x, s\}, b)$-paths internally disjoint from $V(H)$ 
and satisfying the length condition. 
On the other hand, 
since $E_{G}(T, D) \neq \emptyset$
and $E_{G}(T, B-b) = \emptyset$ by the assumption, 
we can take a $(T, b)$-path in $G[V(D - B) \cup \{b\} \cup T]$. 
Combining the $(T, b)$-path with the above
$k- l+2$ $(S \setminus \{x, s\}, b)$-paths, 
we can get 
$k- l+2$ $(T, S \setminus \{x, s\})$-paths internally disjoint from $V(H)$ 
and satisfying the length condition, 
which contradicts Lemma~\ref{lemma:paths satisfying the length condition (1)}(\ref{T to S-{x,s}}). 
\qed

\begin{subclaim}
\label{subclaim:E_{G}({x, s}, D) = emptyset} 
$E_{G}(\{x, s\}, D) = \emptyset$. 
\end{subclaim}
\proof 
Suppose that $E_{G}(\{x, s\}, D) \neq \emptyset$. 
We define the graph $D^{*}$ as follows: 
\begin{align*}
V(D^{*}) 
&= V(D) \cup \{x^{*}, T\} \textup{ and} \\
E(D^{*})
&= E(D) \cup \{vx^{*} : v \in V(D), E_{G}(v, \{x, s\}) \neq \emptyset\} 
\cup \{vT : v \in V(D), E_{G}(v, T) \neq \emptyset\}. 
\end{align*}
Then by Subclaim~\ref{subclaim:E_{G}(B-b, T cup x, s)} and since $E_{G}(T, D) \neq \emptyset$, 
it follows that $(D^{*}, x^{*}, T)$ is a $2$-connected rooted graph. 
By (C3), and since $e_{G}(v, T) \le l$ for $v \in V(D)$ by Claim~\ref{claim:e_{G}(v, T) <= l}(i), 
it also follows that 
each $v \in V(D^{*}) \setminus \{x^{*}, T\}$ satisfies the following: 
\begin{align*}
\deg_{D^{*}}(v) \ge 
\left \{
\begin{array}{ll} 
\deg_{G}(v) - \big( (l-1) + l \big) + 1 \ge 2(k-l + 1)  
& \textup{(if $E_{G}(v, \{x, s\}) = \emptyset$, $E_{G}(v, T) \neq \emptyset$)} \\[3mm]
\deg_{G}(v) - (l + l) + (1 + 1) \ge 2(k-l + 1) 
& \textup{(if $E_{G}(v, \{x, s\}) \neq \emptyset$, $E_{G}(v, T) \neq \emptyset$)}\\[3mm]
\deg_{G}(v) - (l-1) \ge 2(k-l + 1)
& \textup{(if $E_{G}(v, \{x, s\}) = \emptyset$, $E_{G}(v, T) = \emptyset$)}\\[3mm]
 \deg_{G}(v) - l + 1 \ge 2(k-l + 1)
& \textup{(if $E_{G}(v, \{x, s\}) \neq \emptyset$, $E_{G}(v, T) = \emptyset$)}
\end{array}
\right.. 
\end{align*}
Thus $\delta(D^{*}, x^{*}, T) \ge 2(k-l + 1)$. 
By the induction hypothesis, 
$D^{*}$ contains $k-l + 1$ $(T, x^{*})$-paths 
satisfying the length condition. 
Then it follows from the definition of $D^{*}$ that 
$G-V(C)$ contains $k - l + 1$ $(T, \{x, s\})$-paths internally disjoint from $V(H)$ 
and satisfying the length condition, 
which contradicts Lemma~\ref{lemma:paths satisfying the length condition (1)}(\ref{T to {x,s}}). 
\qed

We divide the rest of the proof of Claim~\ref{claim:neighbor of T} 
into two cases as follows.

\medskip

\noindent\textbf{Case~(i) $|N_{G}(D) \cap T| \le 1$. }

\medskip

Since $E_{G}(T, D) \neq \emptyset$, 
we have $|N_{G}(D) \cap T| = 1$, say $N_{G}(D) \cap T = \{t\}$. 
Since $G$ is $2$-connected,
it follows from Subclaim~\ref{subclaim:E_{G}({x, s}, D) = emptyset}
that $E_{G}(D, S \setminus \{x, s\}) \neq \emptyset$.
In particular, $S \setminus \{x, s\} \neq \emptyset$, that is, $l \ge 2$. 
We define the graph $D^{*}$ as follows: 
\begin{align*}
V(D^{*}) 
&= V(D) \cup \{t, S^{*}\} \textup{ and} \\
E(D^{*})
&= E(D) \cup \{vt : v \in V(D), vt \in E(G)\} \cup \{vS^{*} : v \in V(D), E_{G}(v, S \setminus \{x, s\}) \neq \emptyset\}. 
\end{align*}
Then 
by Subclaims~\ref{subclaim:E_{G}(B-b, T cup x, s)} and \ref{subclaim:E_{G}({x, s}, D) = emptyset}, 
$(D^{*}, t, S^{*})$ is a $2$-connected rooted graph. 
By Subclaim~\ref{subclaim:E_{G}({x, s}, D) = emptyset}
and since $l \ge 2$, it also follows that 
for each $v \in V(D^{*}) \setminus \{t, S^{*}\}$, 
\begin{align*}
\deg_{D^{*}}(v) \ge 
\left \{
\begin{array}{ll} 
\deg_{G}(v) \ge 2k  & \textup{(if $E_{G}(v, S \setminus \{x, s \}) = \emptyset$)} \\[3mm]
\deg_{G}(v) - (l-1) + 1 \ge 2 \big( k - l+ 2 ) & \textup{(if $E_{G}(v, S \setminus \{x, s\}) \neq \emptyset$)}
\end{array}
\right.,  
\end{align*}
and thus $\delta(D^{*}, t, S^{*}) \ge 2( k - l+2)$. 
By the induction hypothesis, 
$D^{*}$ contains $k - l+2$ $(t, S^{*})$-paths satisfying the length condition. 
Then it follows from the definition of $D^{*}$ that 
$G-V(C)$ contains $k - l+2$ $(t, S \setminus \{x, s\})$-paths internally disjoint from $V(H)$ 
and satisfying the length condition, 
which contradicts Lemma~\ref{lemma:paths satisfying the length condition (1)}(\ref{T to S-{x,s}}). 

\medskip

\noindent\textbf{Case~(ii) $|N_{G}(D) \cap T| \ge 2$.}

\medskip

Fix $t \in N_{G}(D) \cap T$. 
Note that 
$E_{G}(D, T \setminus \{t\}) \neq \emptyset$. 
We define the graph $D^{*}$ as follows: 
\begin{align*}
V(D^{*}) 
&= V(D) \cup \{t, T^{*}\} \textup{ and} \\
E(D^{*})
&= E(D) \cup \{vt : v \in V(D), vt \in E(G)\} \cup \{vT^{*} : v \in V(D), E_{G}(v, T \setminus \{t\}) \neq \emptyset\}. 
\end{align*}
Then 
by Subclaims~\ref{subclaim:E_{G}(B-b, T cup x, s)} and \ref{subclaim:E_{G}({x, s}, D) = emptyset}, 
$(D^{*}, t, T^{*})$ is a $2$-connected rooted graph. 
It also follows from Claim~\ref{claim:e_{G}(v, T) <= l}(i) 
and Subclaim~\ref{subclaim:E_{G}({x, s}, D) = emptyset} that 
for each $v \in V(D^{*}) \setminus \{t, T^{*}\}$, 
\begin{align*}
\deg_{D^{*}}(v) \ge 
\left \{
\begin{array}{ll} 
\deg_{G}(v) - (l-1) \ge 2(k-l+1) & \textup{(if $E_{G}(v, T \setminus \{t\}) = \emptyset$)} \\[3mm]
\deg_{G}(v) - \big( (l-1) + l \big) + 1 \ge 2(k-l + 1) & \textup{(if $E_{G}(v, T \setminus \{t\}) \neq \emptyset$)}
\end{array}
\right.,  
\end{align*}
and thus $\delta(D^{*}, t, T^{*}) \ge 2 (k-l+1)$. 
By the induction hypothesis, 
$D^{*}$ contains $k-l+1$ $(t, T^{*})$-paths satisfying the length condition. 
Then it follows from the definition of $D^{*}$ that 
$G-V(C)$ contains $k - l + 1$ $T$-paths internally disjoint from $V(H)$ 
and satisfying the length condition, 
which contradicts Lemma~\ref{lemma:paths satisfying the length condition (1)}(\ref{T to T}).

This completes the proof of Claim~\ref{claim:neighbor of T}. 
\qed

\begin{claim}
\label{claim:E_{G}(S -x, C) neq emptyset} 
If $E_{G}(T, C - y) = \emptyset$,
then 
$E_{G}(S \setminus \{x\}, C) = \emptyset$. 
\end{claim}
\proof 
Assume that $E_{G}(T, C - y) = \emptyset$ and $E_{G}(S \setminus \{x\}, C) \neq \emptyset$.
Let $t$ be a vertex of $T$, 
and let $\alpha = 1$ if $ty \in E(G)$; 
otherwise, let $\alpha = 0$. 
Then by Claim~\ref{claim:neighbor of T} and since $E_{G}(T, C-y) = \emptyset$, 
it follows that $N_{G}(t) \subseteq V(H) \cup \{y\}$. 
By the definition of $\alpha$ 
and Lemma~\ref{lemma:l >= k or l >= k-1}(ii), 
we have $l \le k - \alpha - 1$. 
Since $N_{G}(t) \cap V(C) \subseteq \{y\}$, 
by 
(C4) and Claim~\ref{claim:e_{G}(v, T) <= l}(ii), 
we also have $e_{G}(t, T \setminus \{t\}) \le l + \alpha$. 
Combining this with (C1), 
it follows from the inequality $l \le k - \alpha - 1$ that 
\begin{align*}
2k \le \deg_{G}(t) 
&= |S| + e_{G}(t, T \setminus \{t\}) + |E(G) \cap \{ty\}| \le (l + 1) + (l + \alpha) + \alpha\\ 
&\le (k - \alpha)+(k - 1) + \alpha  = 2k - 1, 
\end{align*}
a contradiction. 
Thus the claim holds.
\qed

We divide the proof of Case~2 into two cases according as $|C| = 1$ or $|C| \ge 2$.

\begin{enumerate}[{\textup{{\bf Case~2.\arabic{enumi}.}}}]
\setcounter{enumi}{0}
\item 
$|C| = 1$, i.e., $V(C) = \{y\}$.  
\end{enumerate}

\begin{claim}
\label{claim:neighbor of x and y} 
$N_{G}(x) = N_{G}(y) = T$. 
\end{claim}
\proof 
Since $V(C) \setminus \{ y \}=\emptyset$, 
it follows from Claim~\ref{claim:E_{G}(S -x, C) neq emptyset}
that $E_{G}(y, S \setminus \{x\}) = \emptyset$. 
Since $xy \notin E(G)$ by Claim~\ref{claim:xy is not an edge}, 
we get $E_{G}(y, S) = \emptyset$. 
This together with $\deg_{G}(x) \le \deg_{G}(y)$ 
implies that $N_{G}(x) = N_{G}(y) = T$. 
\qed

\begin{claim}
\label{claim:|T| >= 3} 
$|T| \ge 3$. 
\end{claim}
\proof 
Suppose that $|T| = 2$, say 
$V(T) = \{t_{1}, t_{2}\}$. 
Then by (C1), $|S| = 2$ also holds, say $S \setminus \{x\} = \{s \}$. 
Let $G' = G - \{x, y\}$. 
By Claim~\ref{claim:neighbor of x and y}, 
$G'$ is connected 
and 
$\deg_{G'}(v) = \deg_{G}(v)$ for $v \in V(G') \setminus \{t_{1}, t_{2}\}$. 
If $G'$ is $2$-connected, 
then obviously $(G', t_{1}, t_{2})$ is a $2$-connected rooted graph. 
On the other hand, 
if $G'$ is not $2$-connected, 
then 
since 
$G$ is $2$-connected 
and $N_{G}(x) = N_{G}(y) = \{t_{1}, t_{2}\}$, 
$t_{1}$ and $t_{2}$ are in different end blocks of $G'$; 
since $\{t_{1}, t_{2}\} \subseteq N_{G}(s)$, $s$ is the cut vertex of $G'$ 
contained in both of the two end blocks of $G'$, 
that is, $(G', t_{1}, t_{2})$ is a $2$-connected rooted graph. 
In either case, 
$(G', t_{1}, t_{2})$ is a $2$-connected rooted graph. 
Therefore, 
by the induction hypothesis, 
we can get $k$ $(t_{1}, t_{2})$-paths $\ora{P_{1}}, \dots, \ora{P_{k}}$ in $G'$ satisfying the length condition. 
Then $xt_{1}\ora{P_{1}}t_{2}y, \dots, x t_{1} \ora{P_{k}}t_{2}y$ are $k$ $(x, y)$-paths in $G$ satisfying the length condition, 
a contradiction. 
\qed

Let $s \in S \setminus \{x\}$ and $t \in T$, 
and let $G' = G - \{s, t\}$. 
Then $\deg_{G'}(v) \ge \deg_{G}(v) -2 \ge 2 (k-1)$ for $v \in V(G') \setminus \{x, y\}$. 

We divide the proof of Case~2.1 into three cases according as the connectivity of $G'$.

\begin{enumerate}[{\textup{{\bf Case~2.\arabic{enumi}.\arabic{enumii}. }}}]
\setcounter{enumi}{0}
\setcounter{enumii}{1}
\item
$G'$ is $2$-connected. 
\end{enumerate}

By applying the induction hypothesis to $(G', x, y)$, 
$G'$ contains $k-1$ $(x, y)$-paths satisfying the length condition. 
Let $\ora{P}$ be the longest path in the $k-1$ $(x, y)$-paths, 
and then 
the $k-1$ $(x, y)$-paths 
together with the $(x, y)$-path $x \ora{P} y^{-} s t y$ 
form $k$ $(x, y)$-paths in $G$ satisfying the length condition, 
a contradiction. 

\begin{enumerate}[{\textup{{\bf Case~2.\arabic{enumi}.\arabic{enumii}.}}}] 
\setcounter{enumi}{0}
\setcounter{enumii}{2}
\item 
$G'$ is not connected. 
\end{enumerate}

In this case, 
there exists a component $D$ of $G'$
such that 
$V(D) \cap (V(H) \cup \{y\}) = \emptyset$, 
and let $D' = G[V(D) \cup \{s, t\}]$. 
Since 
$G$ is $2$-connected, 
it follows that 
$(D', s, t)$ is a $2$-connected rooted graph 
such that 
$\deg_{D'}(v) = \deg_{G}(v)$ for $v \in V(D') \setminus \{s, t\} \ (= V(D))$. 
Then by the induction hypothesis, 
$D'$ contains $k$ $(s, t)$-paths satisfying the length condition. 
By adding $xt$ and $st'y$ to each path, where $t'$ is a vertex of $T \setminus \{t\}$, 
we can obtain $k$ $(x, y)$-paths in $G$ satisfying the length condition, 
a contradiction. 

\begin{enumerate}[{\textup{{\bf Case~2.\arabic{enumi}.\arabic{enumii}.}}}] 
\setcounter{enumi}{0}
\setcounter{enumii}{3}
\item 
$G'$ is connected, but not $2$-connected.
\end{enumerate}

By Claims~\ref{claim:neighbor of x and y} 
and \ref{claim:|T| >= 3},
$G[(V(H) \cup \{y\}) \setminus \{s, t\}]$ 
is $2$-connected. 
Since $G'$ is connected but it is not $2$-connected, 
this implies that 
there is an end block $B$ of $G'$ with cut vertex $b$ 
such that $V(B - b) \cap \big( (V(H) \cup \{y\}) \setminus \{s, t\} \big) = \emptyset$. 
Let $\ora{P}$ be a path from $b$ to some vertex $a \in (V(H) \cup \{y\}) \setminus \{s, t\}$ in $G'$ 
internally disjoint from $\big( V(B \cup H) \cup \{y\} \big) \setminus \{s, t\}$. 
Since $N_{G}(x) = N_{G}(y) = T$ by Claim~\ref{claim:neighbor of x and y}, 
it follows that $a \notin \{x, y\}$ 
and thus $a \in V(H) \setminus \{x, s, t\}$. 
Note that $N_{G}(B - b) \subseteq \{s, t, b\}$.

We now show that $t \in N_{G}(B-b)$. 
By way of contradiction, 
suppose $t \notin N_{G}(B-b)$, 
and let 
$B' = G[V(B) \cup \{s\}]$. 
Then $(B', s, b)$ is a $2$-connected rooted graph 
such that 
$\deg_{B'}(v) = \deg_{G}(v)$ for $v \in V(B') \setminus \{s, b\}$. 
Hence by the induction hypothesis, 
$B'$ contains $k$ $(s, b)$-paths $\ora{P_{1}}, \dots, \ora{P_{k}}$ 
satisfying the length condition. 
If $a \in T$, 
then let $\ora{P'} = b \ora{P} a y$; 
if $a \in S$, 
then we take a vertex $t'$ of $T \setminus \{t\}$, 
and let $\ora{P'} = b \ora{P} a t' y$. 
Then 
$x t s \ora{P_{1}} b \ora{P'} y, \dots, x t s \ora{P_{k}} b \ora{P'} y$ 
form $k$ $(x, y)$-paths in $G$ satisfying the length condition, a contradiction. 
Thus $t \in N_{G}(B-b)$ is proved.

Now let $B'' = G[V(B) \cup \{t\}]$. 
Then $(B'', t, b)$ is a $2$-connected rooted graph 
such that 
$\deg_{B''}(v) \ge \deg_{G}(v) - 1 > 2(k-1)$ for $v \in V(B'') \setminus \{t, b\}$. 
By the induction hypothesis, 
$B''$ contains $k-1$ $(t, b)$-paths $\ora{Q_{1}}, \dots, \ora{Q_{k-1}}$ satisfying the length condition. 
If $a \in T$,  
then 
there exists a vertex $t' \in T \setminus \{t, a\}$ 
and thus 
$xt \ora{Q_{1}} b \ora{P} a y, \dots, x t \ora{Q_{k-1}} b \ora{P} a y$ 
and 
$xt \ora{Q_{k-1}} b \ora{P} a s t' y$ 
form $k$ $(x, y)$-paths in $G$ satisfying the length condition, a contradiction. 
Thus $a \in S$. 
Then 
there exist two distinct vertices $t_{1}, t_{2} \in T \setminus \{t\}$, 
and hence 
$x t \ora{Q_{1}} b \ora{P} a t_{1} y, \dots, x t \ora{Q_{k-1}} b \ora{P} a t_{1} y$ 
and 
$x t \ora{Q_{k-1}} b \ora{P} a t_{1} s t_{2} y$ 
form $k$ $(x, y)$-paths in $G$ satisfying the length condition, a contradiction.

This completes the proof of Case~2.1.

\bigskip
\begin{enumerate}[{\textup{{\bf Case~2.\arabic{enumi}.}}}]
\setcounter{enumi}{1}
\item 
$|C| \ge 2$. 
\end{enumerate}

\begin{claim}
\label{claim:E_{G}(T, C-y) neq emptyset} 
$E_{G}(T, C - y) \neq \emptyset$. 
\end{claim}
\proof 
By way of contradiction, suppose that $E_{G}(T, C-y) = \emptyset$. 
By Claim~\ref{claim:E_{G}(S -x, C) neq emptyset},
$E_{G}(S \setminus \{x\}, C) = \emptyset$. 
Hence $N_{G}(C - y) \subseteq \{x, y\}$. 
Let $C' = G[V(C) \cup \{x\}]$. 
Since $G$ is $2$-connected, $(C', x, y)$ is a $2$-connected rooted graph 
such that 
$\deg_{C'}(v) = \deg_{G}(v)$ for $v \in V(C') \setminus \{x, y\}$. 
By the induction hypothesis, 
$C'$ (and also $G$) contains $k$ $(x, y)$-paths satisfying the length condition, 
a contradiction. 
Thus $E_{G}(T, C - y) \neq \emptyset$. 
\qed

By (C1), 
for a vertex $t$ of $T$, 
$H$ contains $2$ $(x, t)$-paths of lengths $1$ and $3$, respectiely. 
Since 
$E_{G}(T, C - y) \neq \emptyset$ by Claim~\ref{claim:E_{G}(T, C-y) neq emptyset}, 
this implies that 
\begin{align}
\label{k >= 3}
k \ge 3. 
\end{align}

In the rest of this proof, 
we say that an end block $B$ of $C$ is \textit{feasible} 
if $y \notin V(B) \setminus \{b\}$, 
where $b$ is the cut vertex of $C$ contained in $B$. 

\begin{claim}
\label{claim:connectivity of C and edges bw T and B-b} 
{\rm (i)} $C$ is not $2$-connected, 
and 
{\rm (ii)} if $B$ is a feasible end block of $C$ 
with cut vertex $b$, 
then $E_{G}(T, B-b) = \emptyset$. 
\end{claim}
\proof 
Suppose that either $C$ is $2$-connected 
or 
there is a feasible end block $B$ of $C$ with cut vertex $b$ 
such that $E_{G}(T, B-b) \neq \emptyset$. 
In the former case, 
we define $B' = C$ and $b' = y$. 
Note that, in this case, $E_{G}(T, B' - b') \neq \emptyset$ holds 
by Claim~\ref{claim:E_{G}(T, C-y) neq emptyset}. 
In the latter case, 
we define $B' = B$ and $b' = b$. 
Note that, in the latter case, $E_{G}(T, B' - b') \neq \emptyset$ also holds 
by the assumption. 
Now we define the graph $B^{*}$ as follows: 
\begin{align*}
V(B^{*}) 
= V(B') \cup \{T\} \textup{ and }
E(B^{*})
= E(B') \cup \{vT : v \in V(B'), E_{G}(v, T) \neq \emptyset\}. 
\end{align*}
Then 
$(B^{*}, T, b')$ is a $2$-connected rooted graph. 
By (C3) and (C4), it also follows that for each $v \in V(B^{*}) \setminus \{T, b'\}$, 
\begin{align*}
\deg_{B^{*}}(v) \ge 
\left \{
\begin{array}{ll} 
\deg_{G}(v)  - l \ge 2(k-l) & \textup{(if $E_{G}(v, T) = \emptyset$)} \\[3mm]
\deg_{G}(v) - \big( l + (l+1) \big) + 1 \ge 2(k-l) & \textup{(if $E_{G}(v, T)  \neq \emptyset$)}
\end{array}
\right.,  
\end{align*}
and thus $\delta(B^{*}, T, b') \ge 2 (k-l)$. 
By the induction hypothesis, 
$B^{*}$ contains $k-l$ $(T, b')$-paths satisfying the length condition. 
Therefore, 
by the definition of $B^{*}$ 
and by adding a $(b', y)$-path in $C$ to each of the $k-l$ paths, 
we can obtain 
$k-l$ $(T, y)$-paths in $G$ 
internally disjoint from $V(H)$ 
and satisfying the length condition, which contradicts Lemma~\ref{lemma:paths satisfying the length condition (2)}(\ref{T to y}). 
\qed

\begin{claim}
\label{claim:E_{G}(S - x, B-b) neq emptyset} 
If $B$ is a feasible end-block of $C$ 
with cut vertex $b$, 
then $E_{G}(S \setminus \{x\}, B-b) \neq \emptyset$.
\end{claim}
\proof 
Suppose that $E_{G}(S \setminus \{x\}, B-b) = \emptyset$. 
By Claim~\ref{claim:connectivity of C and edges bw T and B-b}(ii) 
and the $2$-connectivity of $G$, 
we have $N_{G}(B-b) = \{x, b\}$. 
Let $B' = G[V(B) \cup \{x\}]$. 
Then $(B', x, b)$ is a $2$-connected rooted graph such that 
$\deg_{B'}(v) = \deg_{G}(v)$ for $v \in V(B') \setminus \{x, b\}$. 
By the induction hypothesis, 
$B'$ contains 
$k$ $(x, b)$-paths satisfying the length condition. 
Combining these $k$ paths and a $(b, y)$-path in $C$, 
we can obtain $k$ $(x, y)$-paths in $G$ satisfying the length condition, 
a contradiction. 
\qed

\begin{claim}
\label{claim:l=1 and the neighbor of B-b} 
{\rm (i)} $l = 1$, 
and 
{\rm (ii)} if $B$ is a feasible end block of $C$ 
with cut vertex $b$, 
then $N_{G}(B-b) = S \cup \{b\}$. 
\end{claim}
\proof 
By Claim~\ref{claim:connectivity of C and edges bw T and B-b}(i), 
$C$ contains a feasible end block. 
Let $B$ be an arbitrary feasible end block of $C$ 
with cut vertex $b$, 
and we show that 
$l = 1$ and $N_{G}(B-b) = S \cup \{b\}$. 
Suppose that either $l \ge 2$, 
or $l = 1$ but $N_{G}(B-b) \neq S \cup \{b\}$. 
We define the graph $B^{*}$ as follows: 
\begin{align*}
V(B^{*}) 
= V(B) \cup \{S^{*}\} \textup{ and }
E(B^{*})
= E(B) \cup \{vS^{*} : v \in V(B), E_{G}(v, S \setminus \{x\}) \neq \emptyset\}. 
\end{align*}
Then by Claim~\ref{claim:E_{G}(S - x, B-b) neq emptyset}, 
$(B^{*}, S^{*}, b)$ is a $2$-connected rooted graph.

We first consider the case of $l \ge 2$. 
Then 
by Claim~\ref{claim:connectivity of C and edges bw T and B-b}(ii) and (C3), 
and since $l \ge 2$,
it follows that  
for each $v \in V(B^{*}) \setminus \{S^{*}, b\}$, 
\begin{align*}
\deg_{B^{*}}(v) \ge 
\left \{
\begin{array}{ll} 
\deg_{G}(v)   \ge 2( k - l + 1 )  & \textup{(if $E_{G}(v, S) = \emptyset$)} \\[3mm]
\deg_{G}(v)  - 1 \ge  2( k - l + 1 )  & \textup{(if $E_{G}(v, S \setminus \{x\}) = \emptyset$ and $E_{G}(v, x) \not= \emptyset$)} \\[3mm]
\deg_{G}(v) - l + 1 \ge  2( k - l + 1 )  & \textup{(if $E_{G}(v, S \setminus \{x\}) \neq \emptyset$ and $E_{G}(v, x) = \emptyset$)}\\[3mm]
\deg_{G}(v) - l + 1-1 \ge  2( k - l + 1 )  & \textup{(if $E_{G}(v, S \setminus \{x\}) \neq \emptyset$ and $E_{G}(v, x) \not= \emptyset$)}
\end{array}
\right., 
\end{align*}
and thus $\delta(B^{*}, S^{*}, b) \ge { 2( k - l + 1 ) }$. 
By the induction hypothesis, 
$B^{*}$ contains  $k - l + 1$  $(S^{*}, b)$-paths satisfying the length condition. 
Then 
by the definition of $B^{*}$ 
and by adding a $(b, y)$-path in $C$ to each of the $k - l + 1$ paths, 
we can obtain 
$k - l + 1$ $(S \setminus \{x\}, y)$-paths in $G$ 
internally disjoint from $V(H)$ 
and satisfying the length condition, which contradicts Lemma~\ref{lemma:paths satisfying the length condition (2)}(\ref{S - x to y}).

We next consider the case of $l = 1$.  
In this case,
$N_{G}(B-b) = (S \setminus \{x\}) \cup \{b\}$ 
since 
$E_{G}(S \setminus \{x\}, B-b) \neq \emptyset$ (by Claim~\ref{claim:E_{G}(S - x, B-b) neq emptyset}) 
and 
$N_{G}(B-b) \neq S \cup \{b\}$. 
Hence 
it follows that 
$\deg_{B^{*}}(v) = \deg_{G}(v)$ 
for $v \in V(B^{*}) \setminus \{S^{*}, b\}$. 
By the induction hypothesis, 
$B^{*}$ contains $k$ $(S^{*}, b)$-paths satisfying the length condition. 
Therefore, 
by the definition of $B^{*}$, 
we can easily see that 
$G$ contains $k$ $(x, y)$-paths in $G$ satisfying the length condition, 
a contradiction. 
\qed

By Claim~\ref{claim:l=1 and the neighbor of B-b}(i), 
$|S| = 2$, say $S \setminus \{x\} = \{s\}$. 
On the other hand, 
by Claim~\ref{claim:l=1 and the neighbor of B-b}(ii), 
if 
$B$ is a feasible end block of $C$ with cut vertex $b$ 
and 
$v$ is a vertex of $S$, 
then 
$(G[V(B) \cup \{v\}], v, b)$ is a $2$-connected rooted graph 
such that 
$\deg_{G[V(B) \cup \{v\}]}(v') \ge \deg_{G}(v') -1 \ge 2k - 1$ for $v' \in V(B) \setminus \{b\}$. 
This together with 
the induction hypothesis implies that 
\begin{align}
\label{k-1 (v, b)-paths}
\begin{array}{cc}
\hspace{-48pt}\textup{$G[V(B) \cup \{v\}]$ contains $k-1$ $(v, b)$-paths 
satisfying the length condition} \\[1mm]
\hspace{12pt}\textup{for any feasible end block $B$ of $C$ with cut vertex $b$ 
and any vertex $v \in S$.}
\end{array}
\end{align}

Now let 
\begin{align*}
U = N_{G}(T) \cap \big( V(C) \setminus \{y\}\big). 
\end{align*}
Note that by Claim~\ref{claim:E_{G}(T, C-y) neq emptyset}, 
$U \neq \emptyset$. 
Let 
\begin{align*}
\textup{$B_{1}, \dots, B_{h}$ be all the feasible end blocks of $C$ 
with cut vertices $b_{1}, \dots, b_{h}$, respectively} 
\end{align*}
(note that by Claim~\ref{claim:connectivity of C and edges bw T and B-b}(i), 
such blocks exist).  
We further let 
\begin{align*}
C' = C - \bigcup_{1 \le i \le h}(V(B_{i}) \setminus \{b_{i}\}).
\end{align*}
Note that $C'$ is connected 
and that 
by Claim~\ref{claim:l=1 and the neighbor of B-b}(ii), 
$C'$ contains all the vertices of $U \cup \{b_{1}, b_{2}, \dots, b_{h}\}$.

In the rest of this proof, 
let $\ora{P_{1}}$ and $\ora{P_{2}}$ be $2$ $(x, b_{1})$-paths in $G[V(B_{1}) \cup \{x\}]$ 
satisfying the length condition 
(note that by (\ref{k >= 3}) and (\ref{k-1 (v, b)-paths}), such two paths exist).

\begin{claim}
\label{claim:end block containing y} 
There exists an end block $B_{y}$ of $C$ with cut vertex $b_{y}$ 
such that $y \in V(B_{y}) \setminus \{b_{y}\}$. 
\end{claim}
\proof 
Suppose that $B_{1}, \dots, B_{h}$ are 
all the end blocks of $C$. 
Since 
$y \in V(C')$ and $U \subseteq V(C') \setminus \{y\}$, 
and since 
the block-cut tree of $C$ has order at least $3$, 
there exist two vertex-disjoint paths $\ora{P}$ and $\ora{Q}$ in $C'$ 
such that 
$\ora{P}$ is a path from $b_{i}$ to some vertex $u \in U$ 
and 
$\ora{Q}$ is a path from $b_{j}$ to $y$ 
for some $i, j$ with $i \neq j$, 
say $i = 1$ and $j = 2$. 

On the other hand, 
it follows from (\ref{k-1 (v, b)-paths}) that 
$G[V(B_{2}) \cup \{s\}]$ contains $k-1$ $(s, b_{2})$-paths $\ora{Q_{1}}, \dots, \ora{Q_{k-1}}$ 
satisfying the length condition. 
By the definition of $U$, 
we can take a vertex $t$ of $T$ such that $tu \in E(G)$. 
Then 
\begin{align*}
\begin{array}{llcll}
x\ora{P_{1}} b_{1} \ora{P} u t s \ora{Q_{1}} b_{2} \ora{Q} y, 
& x\ora{P_{1}} b_{1} \ora{P} u t s \ora{Q_{2}} b_{2} \ora{Q} y, 
& \dots, 
& x\ora{P_{1}} b_{1} \ora{P} u t s \ora{Q_{k-1}} b_{2} \ora{Q} y, 
& x\ora{P_{2}} b_{1} \ora{P} u t s \ora{Q_{k-1}} b_{2} \ora{Q} y
\end{array}
\end{align*}
are $k$ $(x, y)$-paths in $G$ satisfying the length condition,
a contradiction. 
\qed

Let $B_{y}$ and $b_{y}$ be the same ones as in Claim~\ref{claim:end block containing y}. 
Then 
$B_{1}, \dots, B_{h}$ and $B_{y}$ are all the end blocks of $C$.

\begin{claim}
\label{claim:e_{G}(v, H) <= 2} 
For each $v \in V(C) \setminus \{y\}$, 
either $e_{G}(v, H) \le 2$ or 
$v$ is a cut vertex of $C$ separating $y$ and all feasible end blocks of $C$. 
\end{claim}
\proof 
Suppose that 
there exist a vertex $v$ of $V(C) \setminus \{y\}$ 
such that 
$e_{G}(v, H) \ge 3$ 
and 
a feasible block $B_{i}$, say $i = 1$, 
such that 
$C - v$ has a $(b_{1}, y)$-path $\ora{Q'}$ 
internally disjoint from $B_{1}$ 
(note that $v \in V(C')$, since every vertex of $\bigcup_{1 \le j \le h}(V(B_{j}) \setminus \{b_{j}\})$ is adjacent to at most two vertices of $H$ 
by Claim~\ref{claim:l=1 and the neighbor of B-b}). 
Since $l = 1$ by Claim~\ref{claim:l=1 and the neighbor of B-b}(i), 
(C3) and (C4) ensure that 
$v$ is adjacent to exactly two distinct vertices in $T$, say $t_{1}$ and $t_{2}$. 
By (\ref{k-1 (v, b)-paths}), 
$G[V(B_{1}) \cup \{s\}]$ contains 
$k-1$ $(s, b_{1})$-paths $\ora{Q_{1}}, \dots, \ora{Q_{k-1}}$ 
satisfying the length condition. 
Then 
the $k-1$ paths $x t_{1} s \ora{Q_{i}} b_{1} \ora{Q'} y$ ($1 \le i \le k - 1$) 
together with the path $x t_{2} v t_{1} s \ora{Q_{k-1}} b_{1} \ora{Q'} y$ 
form $k$ $(x, y)$-paths in $G$ satisfying the length condition, a contradiction. 
\qed

By adding a $(b_{1}, b_{y})$-path in $C'$ to each of $\ora{P_{1}}$ and $\ora{P_{2}}$, 
we can get 
two $(x, b_{y})$-paths $\ora{P_{1}'}$ and $\ora{P_{2}'}$ 
in 
$G[(V(C) \cup \{x\}) \setminus V(B_{y} - b_{y})]$ 
satisfying the length condition. 

\begin{claim}
\label{claim:|B_{y}| = 2} 
$|B_{y}| = 2$ $($i.e., $V(B_{y}) = \{y, b_{y}\})$. 
\end{claim}
\proof 
Suppose that $|B_{y}| \ge 3$, that is, $B_{y}$ is $2$-connected. 
For each vertex $v$ of $V(B_{y}) \setminus \{y, b_{y}\}$, 
$v$ is not a cut vertex of $C$ separating $y$ and all feasible end blocks of $C$, 
and hence by Claim~\ref{claim:e_{G}(v, H) <= 2}, 
we have $e_{G}(v, H) \le 2$. 
This implies that 
$(B_{y}, b_{y}, y)$ is a $2$-connected rooted graph 
such that 
$\deg_{B_{y}}(v) \ge \deg_{G}(v) - 2 \ge 2(k-1)$ for $v \in V(B_{y}) \setminus \{y, b_{y} \}$. 
By the induction hypothesis, 
$B_{y}$ contains $k-1$ $(b_{y}, y)$-paths satisfying the length condition. 
Concatenating these $k-1$ paths with $P_{1}'$ and $P_{2}'$, 
we can obtain $k$ $(x, y)$-paths in $G$ satisfying the length condition, a contradiction. 
\qed

By Claim~\ref{claim:|B_{y}| = 2} and since 
$G$ is $2$-connected, 
$E_{G}(y, H) \neq \emptyset$. 
Since 
$xy \notin E(G)$, 
there exists a vertex $a$ of $V(H) \setminus \{x\}$ 
such that $ay \in E(G)$.

\begin{claim}
\label{claim:h = 1} 
$h = 1$. 
\end{claim}
\proof 
Suppose that $h \ge 2$. 
By (\ref{k-1 (v, b)-paths}), 
$G[V(B_{2}) \cup \{s\}]$ contains $k-1$ $(b_{2}, s)$-paths $\ora{Q_{1}}, \dots \ora{Q_{k-1}}$ 
satisfying the length condition. 
Let $\ora{R}$ be a $(b_{1}, b_{2})$-path in $C'$. 
Note that by Claim~\ref{claim:end block containing y}, 
$y \notin V(R)$. 
We also let 
$\ora{R'}$ be an $(s, a)$-path in $H - x$. 
Then 
\begin{align*}
\begin{array}{llcll}
x \ora{P_{1}} b_{1} \ora{R} b_{2} \ora{Q_{1}} s \ora{R'} a y, 
& x \ora{P_{1}} b_{1} \ora{R} b_{2} \ora{Q_{2}} s \ora{R'} a y, 
& \dots, 
& x \ora{P_{1}} b_{1} \ora{R} b_{2} \ora{Q_{k-1}} s \ora{R'} a y,
& x \ora{P_{2}} b_{1} \ora{R} b_{2} \ora{Q_{k-1}} s \ora{R'} a y 
\end{array}
\end{align*}
are $k$ $(x, y)$-paths in $G$ satisfying the length condition, a contradiction. 
\qed

By Claim~\ref{claim:h = 1}, 
$B_{1}$ and $B_{y}$ are all the end blocks of $C$. 
Therefore, 
$C$ has a unique block $W$ of $C$ 
such that $W \neq B_{y}$ and $b_{y} \in V(W)$. 
Then 
by (C3), (C4), (\ref{k >= 3}), Claims~\ref{claim:l=1 and the neighbor of B-b}(i) and \ref{claim:|B_{y}| = 2}, 
$\deg_{W}(b_{y}) \ge 2k - (2l+1) - |\{b_{y}y\}| \ge 2$, 
which implies that $W$ is $2$-connected. 

We show that $W \not= B_{1}$.
Assume not.
Then by the definition of $U$, Claims~\ref{claim:l=1 and the neighbor of B-b}(ii),~\ref{claim:|B_{y}| = 2} and \ref{claim:h = 1}, 
we have $U \subseteq \{b_{y}\}$.
Then by the definition of $U$, 
we have 
$N_{G}(t) \cap V(C) \subseteq \{y, b_{y}\}$ for $t \in T$. 
Let $t$ be an arbitrary vertex of $T$. 
Since $E_{G}(S \setminus \{x\}, C) \neq \emptyset$ 
by Claim~\ref{claim:l=1 and the neighbor of B-b}(ii), 
it follows from Claim~\ref{claim:neighbor of T} that 
$N_{G}(t) \subseteq V(H \cup C)$. 
Combining this with the above, 
we have $N_{G}(t) \subseteq V(H) \cup \{y, b_{y}\}$. 
Then by (C1), (C4), (\ref{k >= 3}) and Claim~\ref{claim:l=1 and the neighbor of B-b}(i), 
\begin{align*}
6 \le 2k \le \deg_{G}(t) 
= |S| + e_{G}(t, T \setminus \{t\}) + |E(G) \cap \{ty, tb_{y}\}| \le 2l + 4 = 6. 
\end{align*}
Thus the equality holds. 
This yields that 
$k = 3$ and also that $e_{G}(t, T \setminus \{t\}) = l+ 1 = 2$ and $ty, tb_{y} \in E(G)$. 
Since $t$ is an arbitrary vertex of $T$, it follows that 
$G[T]$ contains an edge $t_{1}t_{2}$, 
$t_{1}b_{y} \in E(G)$ and $t_{2}y \in E(G)$. 
Then $x \ora{P_{1}'} b_{y} y , x \ora{P_{2}'} b_{y} y$ and $x \ora{P_{2}'} b_{y} t_{1} t_{2} y$ 
are $3 \ (= k)$ $(x, y)$-paths in $G$ satisfying the length condition, a contradiction. 
Thus 
$W \neq B_{1}$. 
In short, 
\begin{align*}
\textup{$W$ is $2$-connected and $W \neq B_{1}$.}
\end{align*}

Let $w$ be a cut vertex of $C$ which is contained in $W$ such that $w \neq b_{y}$. 
Note that $w$ and $b_{y}$ are all the cut vertices of $C$ which is contained in $W$. 
Then 
each vertex $v$ of $V(W) \setminus \{w, b_{y}\}$ 
is not a cut vertex of $C$ separating $y$ and all feasible end blocks of $C$, 
since $W$ is $2$-connected. 
Hence by Claim~\ref{claim:e_{G}(v, H) <= 2}, 
we have $e_{G}(v, H) \le 2$ for $v \in V(W) \setminus \{w, b_{y}\}$. 
This implies that 
$(W, w, b_{y})$ is a $2$-connected rooted graph 
such that 
$\deg_{W}(v) \ge \deg_{G}(v) - 2 \ge 2(k-1)$ for $v \in V(W) \setminus \{w, b_{y}\}$. 
By the induction hypothesis, 
$W$ contains $k-1$ $(w, b_{y})$-paths satisfying the length condition. 
Concatenating these $k-1$ $(w, b_{y})$-paths with $P_{1}$ and $P_{2}$, the edge $b_{y}y$ 
and a $(b_{1}, w)$-path in $C$, 
we can obtain $k$ $(x, y)$-paths in $G$ satisfying the length condition, a contradiction. 

This completes the proof of Theorem~\ref{improvement of LM2018 Lemma 3.1}. 
\qed



\section{Proof of Theorem~\ref{improvement of LM2018 Lemma 3.1 No.2}}
\label{sec:proof of improvement of LM2018 Lemma 3.1 No.2}

In this section, we give the proof of Theorem~\ref{improvement of LM2018 Lemma 3.1 No.2}. 
The direction is the same as the proof of Theorem~\ref{improvement of LM2018 Lemma 3.1} 
and the argument is also similar. 
Therefore we mainly describe the difference from the proof of Theorem~\ref{improvement of LM2018 Lemma 3.1}. 
In the following proof of Theorem~\ref{improvement of LM2018 Lemma 3.1 No.2}, 
the claims without proof are obtained by the same arguments as in the proof of Theorem~\ref{improvement of LM2018 Lemma 3.1} 
(note that the numberings of the claims correspond to the ones of the claims in the proof of Theorem~\ref{improvement of LM2018 Lemma 3.1}). 

\medskip
\noindent
\textit{Proof of Theorem~\ref{improvement of LM2018 Lemma 3.1 No.2}.}~We prove it by induction on $|V(G)| + |E(G)|$. 
Let $(G, x, y)$ be a minimum counterexample with respect to $|V(G)| + |E(G)|$. 
If $k = 1$, then by (R1) and (R2), 
we can easily see that $G$ contains an $(x, y)$-path of length at least $2$, 
a contradiction. 
Thus $k \ge 2$. 
Since $\delta(G, x, y) \ge {2k-1}$, 
this implies that $|G| \ge 4$. 
By symmetry, we may assume that 
$\deg_{G}(x) \le \deg_{G}(y)$.


\begin{claim}
\label{claim:2-connected No.2} 
$G$ is $2$-connected. 
\end{claim}

\begin{claim}
\label{claim:xy is not an edge No.2} 
$xy \notin E(G)$. 
\end{claim}

\begin{enumerate}[{\textup{{\bf Case~\arabic{enumi}.}}}]
\setcounter{enumi}{0}
\item 
$G-y$ does not contain a cycle of length $4$ passing through $x$. 
\end{enumerate}

Since $xy \notin E(G)$ by Claim~\ref{claim:xy is not an edge No.2}, 
in this case, 
we have 
\begin{align}
\label{e(v, N_{G}(x)) No.2}
e_{G}(v, N_{G}(x) \setminus \{v\}) \le 1 \textup{ for } v \in V(G) \setminus \{x, y\}. 
\end{align}
Let $G^{*}$ be the graph obtained from $G$ 
by contracting the subgraph induced by $N_{G}(x) \cup \{x\}$ 
into a single vertex $x^{*}$ and then removing multiple edges. 
Then by (\ref{e(v, N_{G}(x))  No.2}), 
\begin{align}
\label{minimum degree of G^{*}  No.2}
\deg_{G^{*}}(v) = \deg_{G}(v) \textup{ for } v \in V(G^{*}) \setminus \{x^{*}, y\}. 
\end{align}

Assume for the moment that 
$|G^{*}| = 2$. 
This implies that 
$V(G) = \{x, y\} \cup N_{G}(x)$. 
Since 
$e_{G}(v, N_{G}(x) \setminus \{v\}) \le 1$ for $v \in N_{G}(x)$, 
each vertex $v$ of $N_{G}(x)$ satisfies 
$3 \le {2k-1} \le \deg_{G}(v) \le e_{G}(v, N_{G}(x) \setminus \{v\}) + |\{xv\}| + |E(G) \cap \{yv\}| \le 3$. 
Thus the equality holds. 
This yields that $k = 2$ 
and also that 
$G[N_{G}(x)]$ contains an edge $v_{1}v_{2}$ 
and 
$v_{1}y, v_{2}y \in (G)$. 
Therefore, $xv_{1}y$ and $xv_{1}v_{2}y$ 
are $2 \ (= k)$ $(x, y)$-paths satisfying the semi-length condition, a contradiction. 
Thus $|G^{*}| \ge 3$.

By (\ref{e(v, N_{G}(x)) No.2}), (\ref{minimum degree of G^{*}  No.2}) 
and since $|G^{*}| \ge 3$, 
we can prove the rest of Case~1 
by the same way as in the paragraphs following (\ref{minimum degree of G^{*}}) 
in Case~1 of the proof of Theorem~\ref{improvement of LM2018 Lemma 3.1}.


\begin{enumerate}[{\textup{{\bf Case~\arabic{enumi}.}}}]
\setcounter{enumi}{1}
\item 
$G-y$ contains a cycle of length $4$ passing through $x$. 
\end{enumerate}

By the assumption of Case~2, 
$G$ contains a bipartite subgraph $H = G[S, T]$ 
such that 
$H$ is complete bipartite, $|T| \ge |S| =: l + 1 \ge 2$, 
$x \in S$, $y \notin V(H)$ (i.e., $H$ satisfies (C1) and (C2)). 
Let $C$ be the component of $G-V(H)$ 
such that $y \in V(C)$. 
Choose $H$ so that 
\begin{enumerate}[{\upshape(a)}]
\item $|S|$ is maximum, 
\item $|T|$ is maximal, subject to (a), 
\item $|C|$ is maximum, subject to (a) and (b), and 
\item $|N_{G}(C) \cap S|$ is minimum, subject to (a), (b) and (c). 
\end{enumerate}




\begin{claim}
\label{claim:l-core No.2} 
$H$ is an $l$-core of $G$ with respect to $(x, y)$. 
\end{claim}


\begin{claim}
\label{claim:e_{G}(v, T) <= l No.2} 
If $E_{G}(S \setminus \{x\}, C) \neq \emptyset$, 
then 
{\rm (i)}~$e_{G}(v, T) \le l$ for $v \in V(G) \setminus V(H \cup C)$, 
and 
{\rm (ii)}~$e_{G}(v, T \setminus \{v\}) \le l$  or $E_{G}(v, C) \neq \emptyset$ for $v \in T$. 
\end{claim}

Note that $l \le k - 1$
by Lemma~\ref{lemma:l >= k or l >= k-1} and Claim~\ref{claim:l-core No.2}.

\begin{claim}
\label{claim:neighbor of T No.2} 
If $E_{G}(S \setminus \{x\}, C) \neq \emptyset$, 
then $N_{G}(T) \subseteq V(H \cup C)$. 
\end{claim}

\begin{claim}
\label{claim:E_{G}(S -x, C) neq emptyset No.2} 
If $E_{G}(T, C - y) = \emptyset$,
then 
$E_{G}(S \setminus \{x\}, C) = \emptyset$. 
\end{claim}
\proof 
Assume that $E_{G}(T, C - y) = \emptyset$ and $E_{G}(S \setminus \{x\}, C) \neq \emptyset$. 
Let $t$ be a vertex of $T$, 
and let $\alpha = 1$ if $ty \in E(G)$; 
otherwise, let $\alpha = 0$. 
Then by Claim~\ref{claim:neighbor of T No.2} and since $E_{G}(T, C-y) = \emptyset$, 
it follows that $N_{G}(t) \subseteq V(H) \cup \{y\}$. 
By the definition of $\alpha$ 
and Lemma~\ref{lemma:l >= k or l >= k-1}(ii), 
we have $l \le k - \alpha - 1$. 
Since $N_{G}(t) \cap V(C) \subseteq \{y\}$, 
by 
(C4) and Claim~\ref{claim:e_{G}(v, T) <= l No.2}(ii), 
we also have $e_{G}(t, T \setminus \{t\}) \le l + \alpha$. 
Combining this with (C1), 
it follows from the inequality $l \le k - \alpha - 1$ that 
\begin{align*}
2k-1 \le \deg_{G}(t) 
&= |S| + e_{G}(t, T \setminus \{t\}) + \alpha \le (l + 1) + (l + \alpha) + \alpha\\ 
&\le (k - \alpha)+(k - 1) + \alpha  = 2k - 1. 
\end{align*}
Thus the equality holds in the above inequality.
This implies that  $l = k - \alpha - 1$ and $e_{G}(t, T \setminus \{t\}) = l + \alpha$. 
Since 
$e_{G}(t, T \setminus \{t\}) = l + \alpha \ge l \ge 1$, 
we have $E(G[T]) \not= \emptyset$.
If $\alpha = 0$,
then $l=k-1$,
and so
Lemma~\ref{lemma:l >= k-1, E(T) is not empty}
implies that
$G$ contains $k$ $(x, y)$-paths satisfying 
{the semi-length condition}, a contradiction. 
Therefore $\alpha = 1$, i.e., $ty \in E(G)$. 
Then 
the equality $e_{G}(t, T \setminus \{t\}) = l + \alpha = l + 1 = k-1$ 
gives $|T| \ge k$. 
Now, since $t$ is an arbitrary vertex of $T$, 
these arguments imply the following: 
\begin{align*}
\textup{$G[T]$ contains an edge $t_{1}t_{2}$ such that $t_{1}y, t_{2}y \in E(G)$ and, $|S| = k - 1$ and $|T| \ge k$.}
\end{align*}
Therefore, we can easily find $k$ $(x, y)$-paths in $G[V(H) \cup \{y\}]$ satisfying 
{the semi-length condition}, a contradiction. 
Thus the claim holds.
\qed

We divide the proof of Case~2 into two cases according as $|C| = 1$ or $|C| \ge 2$.

\begin{enumerate}[{\textup{{\bf Case~2.\arabic{enumi}.}}}]
\setcounter{enumi}{0}
\item 
$|C| = 1$, i.e., $V(C) = \{y\}$.  
\end{enumerate}


\begin{claim}
\label{claim:neighbor of x and y No.2} 
$N_{G}(x) = N_{G}(y) = T$. 
\end{claim}

\begin{claim}
\label{claim:|T| >= 3 No.2} 
$|T| \ge 3$. 
\end{claim}

Note that if 
$P_{1}, \dots, P_{k-1}$ are $k-1$ paths satisfying the semi-length condition 
and 
$Q$ is another path of length $|E(P_{k-1})| + 2$, 
then 
$P_{1}, \dots, P_{k-1}, Q$ 
are $k$ paths satisfying the semi-length condition. 
Therefore 
we can prove the rest of Case~2.1 
by the same way as in the paragraphs following Claim~\ref{claim:|T| >= 3} 
in Case~2.1 of the proof of Theorem~\ref{improvement of LM2018 Lemma 3.1}.

\bigskip
\begin{enumerate}[{\textup{{\bf Case~2.\arabic{enumi}.}}}]
\setcounter{enumi}{1}
\item 
$|C| \ge 2$. 
\end{enumerate}


\begin{claim}
\label{claim:E_{G}(T, C-y) neq emptyset No.2} 
$E_{G}(T, C - y) \neq \emptyset$. 
\end{claim}

By (C1), 
for a vertex $t$ of $T$, 
$H$ contains $2$ $(x, t)$-paths of lengths $1$ and $3$, respectively. 
Since 
$E_{G}(T, C - y) \neq \emptyset$ by Claim~\ref{claim:E_{G}(T, C-y) neq emptyset No.2}, 
this implies that 
\begin{align}
\label{k >= 3 No.2}
k \ge 3. 
\end{align}

In the rest of this proof, 
we say that an end-block $B$ of $C$ is \textit{feasible} 
if $y \notin V(B) \setminus \{b\}$, 
where $b$ is the cut vertex of $C$ contained in $B$. 


\begin{claim}
\label{claim:connectivity of C and edges bw T and B-b No.2} 
{\rm (i)} $C$ is not $2$-connected, 
and 
{\rm (ii)} if $B$ is a feasible end block of $C$ 
with cut vertex $b$, 
then $E_{G}(T, B-b) = \emptyset$. 
\end{claim}

\begin{claim}
\label{claim:E_{G}(S - x, B-b) neq emptyset No.2} 
If $B$ is a feasible end block of $C$ 
with cut vertex $b$, 
then $E_{G}(S \setminus \{x\}, B-b) \neq \emptyset$. 
\end{claim}

\begin{claim}
\label{claim:l=1 and the neighbor of B-b No.2} 
{\rm (i)} $l = 1$, 
and 
{\rm (ii)} if $B$ is a feasible end block of $C$ 
with cut vertex $b$, 
then $N_{G}(B-b) = S \cup \{b\}$. 
\end{claim}

By Claim~\ref{claim:l=1 and the neighbor of B-b No.2}(i), 
$|S| = 2$, say $S \setminus \{x\} = \{s\}$. 
On the other hand, 
by Claim~\ref{claim:l=1 and the neighbor of B-b No.2}(ii), 
if 
$B$ is a feasible end block of $C$ with cut vertex $b$ 
and 
$v$ is a vertex of $S$, 
then 
$(G[V(B) \cup \{v\}], v, b)$ is a $2$-connected rooted graph 
such that 
$\deg_{G[V(B) \cup \{v\}]}(v') \ge \deg_{G}(v') -1 \ge 2(k-1)$ for $v' \in V(B) \setminus \{b\}$. 
Therefore, by Theorem~\ref{improvement of LM2018 Lemma 3.1}, 
\begin{align}
\label{k-1 (v, b)-paths No.2}
\begin{array}{cc}
\hspace{-48pt}\textup{$G[V(B) \cup \{v\}]$ contains $k-1$ $(v, b)$-paths 
satisfying {the length condition}} \\[1mm]
\hspace{12pt}\textup{for any feasible end block $B$ of $C$ with cut vertex $b$ 
and any vertex $v \in S$.}
\end{array}
\end{align}
Notice that
these paths satisfy
{the length condition},
not
{the length condition} or {the semi-length condition}. 
Now let 
\begin{align*}
U = N_{G}(T) \cap \big( V(C) \setminus \{y\}\big). 
\end{align*}
Note that by Claim~\ref{claim:E_{G}(T, C-y) neq emptyset No.2}, 
$U \neq \emptyset$. 
Let 
\begin{align*}
\textup{$B_{1}, \dots, B_{h}$ be all the feasible end blocks of $C$ 
with cut vertices $b_{1}, \dots, b_{h}$, respectively} 
\end{align*}
(note that by Claim~\ref{claim:connectivity of C and edges bw T and B-b No.2}(i), 
such blocks exist).  
We further let 
\begin{align*}
C' = C - \bigcup_{1 \le i \le h}(V(B_{i}) \setminus \{b_{i}\}).
\end{align*}
Note that $C'$ is connected 
and that 
by Claim~\ref{claim:l=1 and the neighbor of B-b No.2}(ii), 
$C'$ contains all the vertices of $U \cup \{b_{1}, b_{2}, \dots, b_{h}\}$.

In the rest of this proof, 
let $\ora{P_{1}}, \ldots, \ora{P_{k-1}}$ be $k-1 \ (\ge 2)$ $(x, b_{1})$-paths
in $G[V(B_{1}) \cup \{x\}]$ 
satisfying the length condition 
(note that by (\ref{k >= 3 No.2}) and (\ref{k-1 (v, b)-paths No.2}), such two paths exist).


\begin{claim}
\label{claim:end block containing y No.2} 
There exists an end block $B_{y}$ of $C$ with cut vertex $b_{y}$ 
such that $y \in V(B_{y}) \setminus \{b_{y}\}$. 
\end{claim}

Let $B_{y}$ and $b_{y}$ be the same ones as in Claim~\ref{claim:end block containing y No.2}. 
Then 
$B_{1}, \dots, B_{h}$ and $B_{y}$ are all the end blocks of $C$.

\begin{claim}
\label{claim:e_{G}(v, H) <= 2 No.2} 
For each $v \in V(C) \setminus \{y\}$, 
either $e_{G}(v, H) \le 2$ or 
$v$ is a cut vertex of $C$ separating $y$ and all feasible end blocks of $C$. 
\end{claim}

By adding a $(b_{1}, b_{y})$-path in $C'$ to each $\ora{P_{i}}$, 
we can get $k-1 \ (\ge 2)$ 
$(x, b_{y})$-paths $\ora{P_{1}'}, \dots, \ora{P_{k-1}'}$ 
in $G[(V(C) \cup \{x\}) \setminus V(B_{y} - b_{y})]$ satisfying the length condition.
Note that if 
$\ora{Q_{1}}, \dots, \ora{Q_{k-1}}$ are $k-1$ $(b_{y}, y)$-paths in $B_{y}$ satisfying the semi-length condition, 
then 
$$x\ora{P_{1}'}b_{y}\ora{Q_{1}}y, \ x\ora{P_{1}'}b_{y}\ora{Q_{2}}y, \ \dots, \ x\ora{P_{1}'}b_{y}\ora{Q_{k-1}}y,
 \ x\ora{P_{2}'}b_{y}\ora{Q_{k-1}}y$$
are $k$ $(x, y)$-paths satisfying the semi-length condition. 
Therefore, 
the following claim is also obtained by the same argument as in the proof of Claim~\ref{claim:|B_{y}| = 2}.

\begin{claim}
\label{claim:|B_{y}| = 2 No.2} 
$|B_{y}| = 2$ $($i.e., $V(B_{y}) = \{y, b_{y}\})$. 
\end{claim}

By Claim~\ref{claim:|B_{y}| = 2 No.2} and since 
$G$ is $2$-connected, 
$E_{G}(y, H) \neq \emptyset$. 
Since 
$xy \notin E(G)$, 
\begin{align*}
\textup{there exists a vertex $a$ of $V(H) \setminus \{x\}$ 
such that $ay \in E(G)$.} 
\end{align*}
Furthermore, 
if $N_{G}(y) \cap N_{G}(b_{y}) \cap V(H - x) \neq \emptyset$, 
say $v \in N_{G}(y) \cap N_{G}(b_{y}) \cap V(H - x)$, 
then $x\ora{P_{1}'}b_{y}y, \dots, x\ora{P_{k-1}'}b_{y}y$ and $x\ora{P_{k-1}'}b_{y} v y$ 
are $k$ $(x, y)$-paths satisfying the semi-length condition, a contradiction. 
Thus we also have 
\begin{align}
\label{no common neighbor of y and b_{y} in H}
N_{G}(y) \cap N_{G}(b_{y}) \cap V(H - x) = \emptyset. 
\end{align}


\begin{claim}
\label{claim:h = 1 No.2} 
$h = 1$. 
\end{claim}

By Claim~\ref{claim:h = 1 No.2}, 
$B_{1}$ and $B_{y}$ are all the end blocks of $C$. 
Therefore, 
$C$ has a unique block $W$ of $C$ 
such that $W \neq B_{y}$ and $b_{y} \in V(W)$.

\begin{claim}
\label{claim:W neq B_1 No.2} 
$W \neq B_{1}$.
\end{claim}
\proof 
Suppose that $W = B_{1}$.
Then by the definition of $U$, 
Claims~\ref{claim:l=1 and the neighbor of B-b No.2}(ii),~\ref{claim:|B_{y}| = 2 No.2} and \ref{claim:h = 1 No.2}, 
we have $U \subseteq \{b_{y}\}$.
By the definition of $U$, 
we have 
$N_{G}(t) \cap V(C) \subseteq \{y, b_{y}\}$ for $t \in T$. 
Let $t$ be an arbitrary vertex of $T$. 
Since $E_{G}(S \setminus \{x\}, C) \neq \emptyset$ 
by Claim~\ref{claim:l=1 and the neighbor of B-b No.2}(ii), 
it follows from Claim~\ref{claim:neighbor of T No.2} that 
$N_{G}(t) \subseteq V(H \cup C)$. 
Combining this with the above, 
we have $N_{G}(t) \subseteq V(H) \cup \{y, b_{y}\}$. 
Moreover, 
by (\ref{no common neighbor of y and b_{y} in H}), 
we also have $|E(G) \cap \{ty, tb_{y}\}| \le 1$. 
Therefore, 
by (C1), (C4), Claim~\ref{claim:l=1 and the neighbor of B-b No.2}(i) and (\ref{k >= 3 No.2}), 
\begin{align*}
5 
\le 2k-1 
&\le \deg_{G}(t)\\ 
&= |S| + e_{G}(t, T \setminus \{t\}) + |E(G) \cap \{ty, tb_{y}\}| 
\le (l + 1) + (l + 1) + 1 = 2l + 3 = 5. 
\end{align*}
Thus the equality holds. 
This yields that 
$k = 3$ and also that $e_{G}(t, T \setminus \{t\}) = l+ 1 = 2$ and $|E(G) \cap \{ty, tb_{y}\}| = 1$. 
Since $t$ is an arbitrary vertex of $T$, 
we have 
\begin{align}
\label{e_{G}(t, T - {t}) = 2}
e_{G}(t, T \setminus \{t\}) = 2 \textup{ and } |E(G) \cap \{ty, tb_{y}\}| = 1 \textup{ for } t \in T. 
\end{align}

Assume first that $N_{G}(y) \cap T \not=\emptyset$. 
We may assume that $a \in N_{G}(y) \cap T$ 
(see the paragraph following Claim~\ref{claim:|B_{y}| = 2 No.2}). 
Let $t_{1}, t_{2} \in N_{G}(a) \cap T$ with $t_{1} \neq t_{2}$ 
(note that by (\ref{e_{G}(t, T - {t}) = 2}), such two vertices exist). 
If $t_{i} b_{y} \in E(G)$ for some $i$ with $i \in \{1, 2\}$, 
then 
$x\ora{P_{1}'}b_{y}y, x\ora{P_{2}'}b_{y}y$ and $x\ora{P_{2}'}b_{y} t_{i} a y$ 
are $3 \ ( = k)$ $(x, y)$-paths satisfying {the length condition}, a contradiction. 
Thus $t_{i} b_{y} \notin E(G)$ for $i \in \{1, 2\}$, 
and hence by (\ref{e_{G}(t, T - {t}) = 2}), 
$t_{i} y \in E(G)$ for $i \in \{1, 2\}$. 
Then 
$xay, xat_{1}y$ and 
$xat_{1}st_{2}y$ are 
$3 \ ( = k)$ $(x, y)$-paths satisfying {the semi-length condition}, a contradiction.
Assume next that $N_{G}(y) \cap T=\emptyset$. 
Then $a \in S$, i.e., $a = s$ since $a \in V(H) \setminus \{x\}$.
Since $N_{G}(y) \cap T=\emptyset$,
$tb_{y} \in E(G)$ for $t \in T$.
Hence 
by taking a vertex $t$ of $T$, 
it follows that 
$x\ora{P_{1}'}b_{y}y, x\ora{P_{2}'}b_{y}y$ and $x\ora{P_{2}'}b_{y} t s y$ 
are $3 \ ( = k)$ $(x, y)$-paths satisfying {the length condition}, a contradiction.  
\qed

Moreover, 
we can show that the following holds. 

\begin{claim}
\label{claim:W is 2-connected} 
$W$ is $2$-connected. 
\end{claim}
\proof 
It suffices to show that $\deg_{W}(b_{y}) \ge 2$. 
By way of contradiction, 
suppose that $\deg_{W}(b_{y}) \le 1$. 
Then 
by (C3), (C4), (\ref{k >= 3 No.2}) and Claim~\ref{claim:l=1 and the neighbor of B-b No.2}(i), 
\begin{align*}
5 \le {2k-1} 
&\le \deg_{G}(b_{y}) \\
&= e_{G}(b_{y}, S) + e_{G}(b_{y}, T) + \deg_{W}(b_{y}) + |\{b_{y}y\}| = l + (l + 1) + 1 + 1 = 2l + 3 = 5. 
\end{align*}
Thus the equality holds. 
This yields that 
$k = 3$ and also that $e_{G}(b_{y}, S) = 1$ and $e_{G}(b_{y}, T) = 2$.

Let $t_{1}, t_{2} \in N_{G}(b_{y}) \cap T$ with $t_{1} \neq t_{2}$. 
If $b_{y}x \in E(G)$, 
then $xb_{y}y, xt_{1}b_{y}y$ and $xt_{1} s t_{2} b_{y} y$ 
are $3 \ ( = k)$ $(x, y)$-paths satisfying {the semi-length condition}, a contradiction.
Thus $b_{y}x \notin E(G)$, 
and hence the equality 
$e_{G}(b_{y}, S) = 1$ implies that $b_{y}s \in E(G)$. 
Then by (\ref{no common neighbor of y and b_{y} in H}), 
$a \in T$, 
and hence 
$x\ora{P_{1}'}b_{y}y, x\ora{P_{2}'}b_{y}y$ 
and $x\ora{P_{2}'}b_{y}say$ are $3 \ ( = k)$ $(x, y)$-paths satisfying {the length condition}, a contradiction. 
Thus $\deg_{W}(b_{y}) \ge 2$. 
\qed

Since $W \neq B_{1}$ and $W$ is $2$-connected by Claims~\ref{claim:W neq B_1 No.2} and \ref{claim:W is 2-connected}, 
we can prove the rest of Case~2.2
by the same way as in the last paragraph of Case~2.2 in the proof of Theorem~\ref{improvement of LM2018 Lemma 3.1}.

This completes the proof of Theorem~\ref{improvement of LM2018 Lemma 3.1 No.2}. 
\qed

\section{Proofs of Theorems~\ref{improvement of LM2018 Lemma 4.1} and \ref{improvement of LM2018 Theorem 5.2}}
\label{sec:proofs of improvement of LM2018 Lemma 4.1 and Theorem 5.2}

In this section, 
for a positive integer $k$, 
we let 
\begin{align*}
\varphi \ (= \varphi(k)) 
= 
\left \{
\begin{array}{ll} 
0 & \textup{(if $k$ is odd)} \\[3mm]
1 & \textup{(if $k$ is even)}
\end{array}
\right..
\end{align*}

Now we first show Theorem~\ref{improvement of LM2018 Lemma 4.1}.

\medskip
\noindent
\textit{Proof of Theorem~\ref{improvement of LM2018 Lemma 4.1}.}~By the definition of $\varphi$, 
we have $k = 2l - 1 + \varphi$ for some $l \ge 1$. 
Then 
by the degree condition, $\delta(G) \ge k + 1 = 2l + \varphi = 2(l + \varphi) - \varphi$. 
Since 
$G$ is $2$-connected but not $3$-connected, 
there exists a separation $(A, B)$ of $G$ 
of order two, say $A \cap B = \{x, y\}$. 
Then it is easily seen that 
each of $(G[A], x, y)$ and $(G[B], x, y)$ 
is a $2$-connected rooted graph and, 
\begin{align*}
\delta(G[A], x, y) \ge \delta(G) \ge 2(l + \varphi) - \varphi  \ge 2(l + \varphi) - 1 
\mbox{ and also } 
\delta(G[B], x, y) \ge 2(l + \varphi) - \varphi. 
\end{align*}
Therefore, 
by applying Theorem~\ref{improvement of LM2018 Lemma 3.1 No.2} 
to 
each of $(G[A], x, y)$ and $(G[B], x, y)$, 
it follows that 
$G[A]$ (resp., $G[B]$) contains $l + \varphi$ $(x, y)$-paths $\ora{P_{1}}, \dots, \ora{P_{l+ \varphi}}$ 
(resp., $\ora{Q_{1}}, \dots, \ora{Q_{l+ \varphi}}$) 
satisfying the length condition or the semi-length condition. 
In particular, 
Theorem~\ref{improvement of LM2018 Lemma 3.1} guarantees that 
\begin{align*}
\textup{both of $P_{1}, \dots, P_{l+ \varphi}$ 
and $Q_{1}, \dots, Q_{l+ \varphi}$ satisfy the length condition if $\varphi = 0$.}
\end{align*}

Now, 
suppose that either $P_{1}, \dots, P_{l+ \varphi}$ or $Q_{1}, \dots, Q_{l+ \varphi}$ satisfy the length condition.
By the symmetry, 
we may assume that 
$P_{1}, \dots, P_{l+ \varphi}$ satisfy the length condition. 
By applying Theorem~\ref{improvement of LM2018 Lemma 3.1} to $(G[B], x, y)$, 
we can take other $l$ $(x, y)$-paths $\ora{Q_{1}'}, \dots, \ora{Q_{l}'}$ 
satisfying the length condition, since $\delta(G[B], x, y) \ge 2(l + \varphi) - \varphi \ge 2l$. 
Then Table~\ref{k cycles for the case 2-connected (1)}~\footnote{Consider in order from the first column of the left in Table~\ref{k cycles for the case 2-connected (1)}}\! 
gives $l+(l +\varphi -1) = 2l - 1 + \varphi =k$ cycles in $G$ satisfying the length condition. 
Therefore, 
we may assume that 
neither $P_{1}, \dots, P_{l+ \varphi}$ nor $Q_{1}, \dots, Q_{l+ \varphi}$ satisfy the length condition, 
and then both of $P_{1}, \dots, P_{l+ \varphi}$ and $Q_{1}, \dots, Q_{l+ \varphi}$ 
satisfy the semi-length condition. 
Note that $\varphi=1$.

Let $p$ and $q$ be the switches of 
$P_{1}, \dots, P_{l+\varphi} \ (= P_{l+ 1})$ and $Q_{1}, \dots, Q_{l+\varphi} \ (= Q_{l+ 1})$, 
respectively. 
Then 
Table~\ref{k cycles for the case 2-connected (2)} 
gives $2l = 2l -1 + \varphi = k$ cycles in $G$ satisfying the length condition. 
\qed

\begin{table}[H]
\begin{center}
\footnotesize
\renewcommand{\arraystretch}{1.5}
\begin{tabularx}{120mm}{ C | C }
\hline
\multicolumn{2}{c}{Both of $P_{1}, \dots, P_{l+ \varphi}$ 
and $Q_{1}', \dots, Q_{l}'$ satisfy the length condition}\\
\hline \hline 
$x \, \ora{P_{1}} \, y \, \ola{Q_{1}'} \, x$ 
&$x \, \ora{P_{2}} \, y \, \ola{Q_{l}'} \, x$
\\
$x \, \ora{P_{1}} \, y \, \ola{Q_{2}'} \, x$ 
&$x \, \ora{P_{3}} \, y \, \ola{Q_{l}'} \, x$
\\
$\vdots$ 
& $\vdots$
\\
$x \, \ora{P_{1}} \, y \, \ola{Q_{l}'} \, x$ 
&$x \, \ora{P_{l+\varphi}} \, y \, \ola{Q_{l}'} \, x$
\\
\hline 
$l$ cycles 
& $l + \varphi - 1$ cycles 
\\
\hline
\end{tabularx}
\caption{}
\label{k cycles for the case 2-connected (1)}
\end{center} 
\end{table}

\begin{table}[H]
\begin{center}
\footnotesize
\renewcommand{\arraystretch}{1.5}
\begin{tabularx}{150mm}{ C | C | C | C | C }
\hline
\multicolumn{5}{c}{Both of $P_{1}, \dots, P_{l+ \varphi}$ 
and $Q_{1}, \dots, Q_{l+ \varphi}$ satisfy the semi-length condition and $\varphi = 1$}\\
\hline \hline 
$x \, \ora{P_{1}} \, y \, \ola{Q_{1}} \, x$ 
&$x \, \ora{P_{2}} \, y \, \ola{Q_{q}} \, x$ 
&  
&$x \, \ora{P_{p+1}} \, y \, \ola{Q_{q + 2}} \, x$
&$x \, \ora{P_{p+2}} \, y \, \ola{Q_{l+1}} \, x$
\\
$x \, \ora{P_{1}} \, y \, \ola{Q_{2}} \, x$ 
&$x \, \ora{P_{3}} \, y \, \ola{Q_{q}} \, x$
&$x \, \ora{P_{p+1}} \, y \, \ola{Q_{q + 1}} \, x$ 
&$\vdots$ 
&$\vdots$
\\
$\vdots$ 
&$\vdots$ 
&  
&$x \, \ora{P_{p+1}} \, y \, \ola{Q_{l}} \, x$ 
&$x \, \ora{P_{l}} \, y \, \ola{Q_{l+1}} \, x$
\\
$x \, \ora{P_{1}} \, y \, \ola{Q_{q}} \, x$ 
&$x \, \ora{P_{p}} \, y \, \ola{Q_{q}} \, x$ 
&  
&$x \, \ora{P_{p+1}} \, y \, \ola{Q_{l+1}} \, x$ 
&$x \, \ora{P_{l+1}} \, y \, \ola{Q_{l+1}} \, x$
\\
\hline
$q$ cycles 
&$p-1$ cycles 
&$1$ cycle 
&$l-q$ cycles 
&$l-p$ cycles 
\\
\hline
\end{tabularx}
\caption{}
\label{k cycles for the case 2-connected (2)}
\end{center} 
\end{table}

We next show Theorem~\ref{improvement of LM2018 Theorem 5.2}. 
In the proof, 
we also use the following two lemmas (Lemmas~\ref{LM2018 Lemma 5.1} and \ref{lemma:consecutive cycles}).

\begin{lemma}[\textup{\cite[Lemma~5.1]{LM2018}}]
\label{LM2018 Lemma 5.1}
Let $G$ be a connected graph such that $\delta(G) \ge 4$. 
If $G$ contains a non-separating induced odd cycle, 
then $G$ contains a non-separating induced odd cycle $\ora{C}$ satisfying 
one of the following {\rm (1)} and {\rm (2)}. 
\begin{enumerate}[{\upshape(1)}]
\item $|V(C)| = 3$, or 
\item 
for every non-cut vertex $v$ of $G-V(C)$, 
$e_{G}(v, C) \le 2$, and the equality holds 
if and only if $vu^{+}, vu^{-} \in E(G)$ for some $u \in V(C)$. 
\end{enumerate}
\end{lemma}

\begin{lemma}
\label{lemma:consecutive cycles} 
Let $k$ and $l$ be positive integers
such that $k = 2l - 1 + \varphi(k)$. 
Let $G$ be a graph and $\ora{C}$ be an odd cycle in $G$, say $|V(C)| = 2m + 1$ for some $m \ge 1$. 
We further let $x \in V(G) \setminus V(C)$ and $u \in V(C)$. 
If one of the following {\rm (i)--(iii)} holds, 
then $G$ contains $k$ cycles having consecutive lengths. 
\begin{enumerate}[{\upshape(i)}]
\item 
\label{u to u^{+m} (1)}
$\varphi(k) = 0$ and $G$ contains $l$ $(u, \{u^{+m}, u^{-m}\})$-paths 
internally disjoint from $V(C)$ 
and satisfying the length condition 
or the semi-length condition. 
\item 
\label{u to u^{+m} (2)}
$\varphi(k) = 1$ and 
$G$ contains $l$ $(u, \{u^{+m}, u^{-m}\})$-paths 
internally disjoint from $V(C)$ 
and satisfying the length condition. 
\item 
\label{x to u^{+m}}
$m\ge2$, 
$xu^{+}, xu^{-} \in E(G)$ 
and 
$G$ contains $l-1$ $(x, u^{+m})$-paths internally disjoint from $V(C) \cup \{x\}$ 
and satisfying the length condition. 
\end{enumerate} 
\end{lemma}

\medskip
\noindent
\textit{Proof of Lemma~\ref{lemma:consecutive cycles}.}~To show (\ref{u to u^{+m} (1)}) and (\ref{u to u^{+m} (2)}), 
suppose that 
$G$ contains $l$ $(u, \{u^{+m}, u^{-m}\})$-paths $\ora{P_{1}}, \dots, \ora{P_{l}}$ 
internally disjoint from $V(C)$ 
and satisfying the length condition 
or the semi-length condition. 
We further suppose that 
if $\varphi = 1$, then $P_{1}, \dots, P_{l}$ 
satisfy the length condition. 
For each $i$ with $1 \le i \le l$, 
let $v_{i}$ be the end vertex of $P_{i}$ such that $v_{i} \in \{u^{+m}, u^{-m}\}$, 
and let $\ora{Q_{i}}$ and $\ora{R_{i}}$ 
denote the paths in $C$ from $v_{i}$ to $u$ 
such that $|E(Q_{i})| = m$ and $|E(R_{i})| = m+1$. 
If 
$P_{1}, \dots, P_{l}$ satisfy the length condition, 
then Table~\ref{k cycles with consecutive lengths (1)} 
gives $2l \ge  2l - 1 + \varphi  = k$ cycles having consecutive lengths (see also Figure~\ref{consecutivecycles1}). 
Thus we may assume that 
$P_{1}, \dots, P_{l}$ does not satisfy the length condition 
but satisfy the semi-length condition. 
In particular, by our assumption, we have $\varphi \neq 1$, i.e., $\varphi = 0$. 
Let $j$ be the switch of $P_{1}, \dots, P_{l}$, 
and then Table~\ref{k cycles with consecutive lengths (2)} gives 
$2l-1 = 2l - 1 + \varphi = k$ cycles having consecutive lengths. 
Thus (\ref{u to u^{+m} (1)}) and (\ref{u to u^{+m} (2)}) are proved.

\begin{figure}[H]
\begin{center}
\input{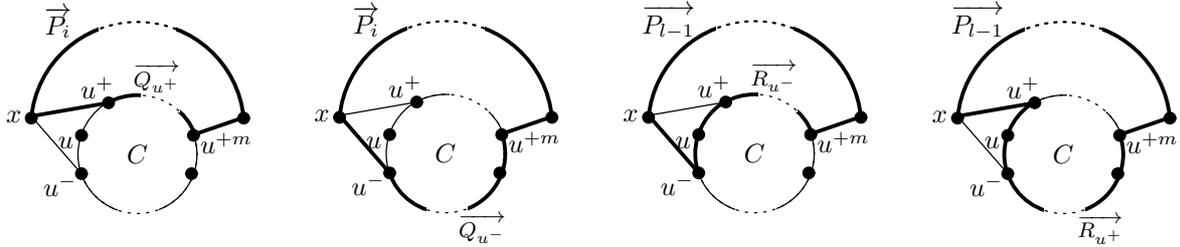}
\caption{The paths $P_{i}$, $Q_{u^+}$, $Q_{u^-}$, $R_{u^+}$ and $R_{u^-}$}
\label{consecutivecycles2}
\end{center}
\end{figure}

We next show (\ref{x to u^{+m}}). 
Suppose that  
$m \ge 2$, $xu^{+}, xu^{-} \in E(G)$ 
and 
$G$ contains $l - 1$ $(x, u^{+m})$-paths 
$\ora{P_{1}}, \dots, \ora{P_{l-1}}$ 
internally disjoint from $V(C) \cup \{x\}$ 
and satisfying the length condition. 
Since $|C| = 2m + 1 \ge 5$,
we can let $\ora{Q_{u^{+}}}$ and $\ora{R_{u^{+}}}$ (resp., $\ora{Q_{u^{-}}}$ and $\ora{R_{u^{-}}}$)
be the paths in $C$ 
from $u^{+m}$ to $u^{+}$ (resp., from $u^{+m}$ to $u^{-}$)
such that 
$|E(Q_{u^{+}})| = m - 1$ and $|E(R_{u^{+}})| = m + 2$ 
(resp., $|E(Q_{u^{-}})| = m$ and $|E(R_{u^{-}})| = m + 1$). 
Hence 
Table~\ref{k cycles with consecutive lengths (3)} 
gives $2l \ge  2l - 1 + \varphi  = k$ cycles having consecutive lengths (see also Figure~\ref{consecutivecycles2}). 
Thus (\ref{x to u^{+m}}) is also proved. 
\qed

\begin{figure}[H]
\begin{minipage}{0.5\hsize}
\begin{center}
\footnotesize
\renewcommand{\arraystretch}{1.5}
\makeatletter
\def\@captype{table}
\makeatother
\begin{tabularx}{60mm}{C}
\hline
$P_{1}, \dots, P_{l}$ satisfy the length condition\\
\hline\hline
$u \, \ora{P_{1}} \, v_{1} \, \ora{Q_{1}} \, u$, \ \ $u \, \ora{P_{1}} \, v_{1} \, \ora{R_{1}} \, u$ \\
$u \, \ora{P_{2}} \, v_{2} \, \ora{Q_{2}} \, u$, \ \ $u \, \ora{P_{2}} \, v_{2} \, \ora{R_{2}} \, u$ \\
$\vdots$ \\
$u \, \ora{P_{l}} \, v_{l} \, \ora{Q_{l}} \, u$, \ \ $u \, \ora{P_{l}} \, v_{l} \, \ora{R_{l}} \, u$ \\
\hline 
$2l$ cycles \\
\hline
\end{tabularx}
\caption{}
\label{k cycles with consecutive lengths (1)}
\end{center}
\end{minipage}
\begin{minipage}{0.5\hsize}
\begin{center}
\input{consecutivecycles1.tex}
\caption{The paths $P_{i}$, $Q_{i}$ and $R_{i}$}
\label{consecutivecycles1}
\end{center}
\end{minipage}
\end{figure}

\begin{table}[H]
\begin{center}
\footnotesize
\renewcommand{\arraystretch}{1.5}
\begin{tabularx}{140mm}{ C | C }
\hline
\multicolumn{2}{c}{$P_{1}, \dots, P_{l}$ 
satisfy the semi-length condition and $\varphi = 0$}\\
\hline \hline 
$u \, \ora{P_{1}} \, v_{1} \, \ora{Q_{1}} \, u$, \ $u \, \ora{P_{1}} \, v_{1} \, \ora{R_{1}} \, u$ 
&$u \, \ora{P_{j+1}} \, v_{j+1}  \, \ora{R_{j+1}} \, u$  
\\
$u \, \ora{P_{2}} \, v_{2} \, \ora{Q_{2}} \, u$, \ $u \, \ora{P_{2}} \, v_{2} \, \ora{R_{2}} \, u$ 
&$u \, \ora{P_{j+2}} \, v_{j+2} \, \ora{Q_{j+2}} \, u$, \ $u \, \ora{P_{j+2}} \, v_{j+2} \, \ora{R_{j+2}} \, u$ 
\\
$\vdots$ 
&$u \, \ora{P_{j+3}} \, v_{j+3} \, \ora{Q_{j+3}} \, u$, \ $u \, \ora{P_{j+3}} \, v_{j+3} \, \ora{R_{j+3}} \, u$ 
\\
$u \, \ora{P_{j-1}} \, v_{j-1} \, \ora{Q_{j-1}} \, u$, \ $u \, \ora{P_{j-1}} \, v_{j-1} \, \ora{R_{j-1}} \, u$ 
&$\vdots$  
\\
$u \, \ora{P_{j}} \, v_{j} \, \ora{Q_{j}} \, u$, \ $u \, \ora{P_{j}} \, v_{j} \, \ora{R_{j}} \, u$ 
&$u \, \ora{P_{l}} \, v_{l} \, \ora{Q_{l}} \, u$, \ $u \, \ora{P_{l}} \, v_{l} \, \ora{R_{l}} \, u$  
\\
\hline 
$2j$ cycles 
& $2l-2j-1$ cycles
\\
\hline
\end{tabularx}
\caption{}
\label{k cycles with consecutive lengths (2)}
\end{center} 
\end{table}

\begin{table}[H]
\begin{center}
\footnotesize
\renewcommand{\arraystretch}{1.5}
\begin{tabularx}{140mm}{ C | C }
\hline
\multicolumn{2}{c}{$m \ge 2$, \ $xu^{+}, xu^{-} \in E(G)$ and $P_{1}, \dots, P_{l-1}$ 
satisfy the semi-length condition}\\
\hline \hline 
$x \, \ora{P_{1}} \, u^{+m} \, \ora{Q_{u^{+}}} \, u^{+} \, x$, \ $x \, \ora{P_{1}} \, u^{+m} \, \ora{Q_{u^{-}}} \, u^{-} \, x$
& 
\\
$x \, \ora{P_{2}} \, u^{+m} \, \ora{Q_{u^{+}}} \, u^{+} \, x$, \ $x \, \ora{P_{2}} \, u^{+m} \, \ora{Q_{u^{-}}} \, u^{-} \, x$
&$x \, \ora{P_{l-1}} \, u^{+m} \, \ora{R_{u^{-}}} \, u^{-} \, x$, \ $x \, \ora{P_{l-1}} \, u^{+m} \, \ora{R_{u^{+}}} \, u^{+} \, x$ 
\\
$\vdots$ 
&
\\
$x \, \ora{P_{l-1}} \, u^{+m} \, \ora{Q_{u^{+}}} \, u^{+} \, x$, \ $x \, \ora{P_{l-1}} \, u^{+m} \, \ora{Q_{u^{-}}} \, u^{-} \, x$
&
\\
\hline 
$2l-2$ cycles 
& $2$ cycles
\\
\hline
\end{tabularx}
\caption{}
\label{k cycles with consecutive lengths (3)}
\end{center} 
\end{table}

We are now ready to prove Theorem~\ref{improvement of LM2018 Theorem 5.2}.

\medskip
\noindent
\textit{Proof of Theorem~\ref{improvement of LM2018 Theorem 5.2}.}~If 
$k = 1$, then the assertion clearly holds, since $G$ is $2$-connected. 
If 
$k = 2$, 
take an edge $xy$ of $G$ 
and 
find $2$ $(x, y)$-paths satisfying the length condition or the semi-length condition 
by using Theorem~\ref{improvement of LM2018 Lemma 3.1 No.2}, 
and then the edge $xy$ and the $2$ $(x, y)$-paths 
induce $2$ cycles which satisfy the length condition or have consecutive lengths. 
Thus we may assume $k \ge 3$.

Let $k = 2l - 1 + \varphi$ for some $l \ge 2$. 
Then by the degree condition, 
we have 
$\delta(G) \ge k + 1 \ (= 2l + \varphi) \ge 4$. 
Hence by Lemma~\ref{LM2018 Lemma 5.1}, 
$G$ contains a non-separating induced odd cycle $\ora{C}$ in $G$, say $|V(C)| = 2m + 1$ for $m \ge 1$, 
such that $C$ satisfies (1) or (2) in Lemma~\ref{LM2018 Lemma 5.1}.
Now, suppose that $G$ does not contain $k$ cycles having consecutive lengths, 
and then 
we will show that $G$ contains $k$ cycles satisfying the length condition.

\begin{claim}
\label{claim:|C| >= 5} 
$|V(C)| \ge 5$, that is, $m \ge 2$. 
\end{claim}
\proof 
Suppose that $|V(C)| = 3$, that is, $m=1$.
Let $u \in V(C)$. 
Consider the graph $G^{*}$ 
obtained from $G$ by contracting $u^{+}$ and $u^{-}$ into a vertex $u^{*}$. 
Then 
$G^{*}$ is a $2$-connected graph and 
$\delta(G^{*}) \ge \delta(G) - 1 \ge k = 2l - 1 + \varphi \ (\ge 2l - 1)$. 
Hence by Theorem~\ref{improvement of LM2018 Lemma 3.1 No.2}, 
$G^{*}$ contains $l$ $(u, u^{*})$-paths 
satisfying 
the length condition or the semi-length condition. 
In particular, if $\varphi = 1$, then by 
Theorem~\ref{improvement of LM2018 Lemma 3.1}, 
we may assume that 
the $l$ $(u, u^{*})$-paths 
satisfy the length condition. 
Note that each of the $l$ paths does not contain the edge $uu^{*}$, 
since the length is at least $2$. 
Then 
it follows from the definition of $G^{*}$ that 
$G$ contains $l$ 
$(u, \{u^{+}, u^{-}\})$-paths internally disjoint from $V(C)$ 
and satisfying the length condition or the semi-length condition; 
in particular,  
they satisfy the former condition when 
$\varphi = 1$. 
This together with 
Lemma~\ref{lemma:consecutive cycles}(\ref{u to u^{+m} (1)}) or (\ref{u to u^{+m} (2)}) leads to a contradiction. 
\qed

By Claim~\ref{claim:|C| >= 5}, 
$C$ satisfies {\rm (2)} in Lemma~\ref{LM2018 Lemma 5.1}, 
i.e., 
every non-cut vertex $v$ of $G-V(C)$ 
satisfies 
\begin{align}\label{deg_{G-V(C)}(v)}
\deg_{G-V(C)}(v) \ge \delta(G) - 2 \ge k - 1 = 2(l-1) + \varphi. 
\end{align}

\begin{fact}\label{B is 2-connected}
If $B$ is an end block of $G-V(C)$,
then $B$ is $2$-connected.
\end{fact}

\proof
Since $k \ge 3$, (\ref{deg_{G-V(C)}(v)}) implies that $|B| \ge 3$, 
and hence  $B$ is $2$-connected. 
\qed

\begin{claim}
\label{claim:neighbor of u^{+m}} 
Let $B$ be an end block of $G-V(C)$ with cut vertex $b$ 
when $G-V(C)$ is not $2$-connected; 
otherwise, let $B = G - V(C)$ 
and $b$ be an arbitrary vertex of $G-V(C)$. 
Further, let $x \in V(B - b)$ and $u \in V(C)$. 
If 
one of the following {\rm (i)} and {\rm (ii)} holds, 
then 
we have 
$E_{G}(\{u^{+m}, u^{-m}\}, V(G-C) \setminus V(B-b)) = \emptyset$. 
\begin{enumerate}[{\upshape(i)}]
\item 
$xu \in E(G)$ 
and 
every vertex of $B-b$ is adjacent to at most one vertex of $C$.  
\item
$xu^{+}, xu^{-} \in E(G)$. 
\end{enumerate}
\end{claim}
\proof 
Suppose that (i) or (ii) holds 
and that $E_{G}(\{u^{+m}, u^{-m}\}, V(G-C) \setminus V(B-b)) \neq \emptyset$. 
By the symmetry of $u^{+m}$ and $u^{-m}$,
we may assume that 
$E_{G}(u^{+m}, V(G-C) \setminus V(B-b)) \neq \emptyset$. 
Let 
$y$ be a vertex of $V(G-C) \setminus V(B-b)$ such that $yu^{+m} \in E(G)$.

We first consider the case where (i) holds. 
Then by the assumption, 
every vertex $v$ of $B - b$ satisfies 
$$\deg_{B}(v) = \deg_{G-V(C)}(v) \ge \delta(G) - 1 \ge k = 2l - 1 + \varphi \ (\ge 2l - 1).$$
Then it follows from Fact~\ref{B is 2-connected} and Theorem~\ref{improvement of LM2018 Lemma 3.1 No.2} 
that 
$B$ contains $l$ $(x, b)$-paths 
satisfying 
the length condition or the semi-length condition. 
In particular, if $\varphi = 1$, then by 
Theorem~\ref{improvement of LM2018 Lemma 3.1}, 
we may assume that 
the $l$ $(x, b)$-paths 
satisfy the length condition. 
By adding a $(b, y)$-path in $G-(V(C) \cup V(B - b))$ 
to each of the $l$ $(x, b)$-paths, 
we can get $l$
$(x, y)$-paths in $G-V(C)$
satisfying the length condition or the semi-length condition; 
in particular, they satisfy the former condition when $\varphi = 1$. 
Since $xu, yu^{+m} \in E(G)$, 
this together with 
Lemma~\ref{lemma:consecutive cycles}(\ref{u to u^{+m} (1)}) or (\ref{u to u^{+m} (2)}) 
leads to a contradiction.

We next consider the case where (ii) holds. 
Then 
it follows from (\ref{deg_{G-V(C)}(v)}), Fact~\ref{B is 2-connected} and Theorem~\ref{improvement of LM2018 Lemma 3.1} 
that 
$B$ contains $l-1$ $(x, b)$-paths 
satisfying the length condition. 
By adding a $(b, y)$-path in $G-(V(C) \cup V(B - b))$ 
to each of the $l-1$ $(x, b)$-paths, 
we can get $l-1$ 
$(x, y)$-paths in $G-V(C)$
satisfying the length condition. 
Since 
$xu^{+}, xu^{-}, yu^{+m} \in E(G)$ 
and since 
$m \ge 2$ by Claim~\ref{claim:|C| >= 5}, 
this together with 
Lemma~\ref{lemma:consecutive cycles}(\ref{x to u^{+m}}) 
leads to a contradiction. 
\qed

\begin{claim}
\label{claim:G - V(C) is not 2-connected} 
$G - V(C)$ is not $2$-connected. 
\end{claim}
\proof 
Suppose that 
$B := G - V(C)$ is $2$-connected. 
Choose a vertex $x$ of $B$ so that 
$e_{G}(x, C)$ is as large as possible. 
Then $e_{G}(x, C) = 1$ or $2$ 
and $e_{G}(v, C) \le 2$ for $v \in V(B)$ 
since $G$ is $2$-connected and $C$ satisfies {\rm (2)} in Lemma~\ref{LM2018 Lemma 5.1}. 
Furthermore, if 
$e_{G}(x, C) = 2$, 
then there is a vertex $u$ of $C$ such that $xu^{+}, xu^{-} \in E(G)$; 
otherwise, 
let $u$ be the unique vertex of $C$ such that $xu \in E(G)$. 
Note that if $e_{G}(x, C) \neq 2$, then the choice of $x$ implies that 
every vertex of $B$ is adjacent to at most one vertex of $C$.  
Now, 
since $\delta(G) \ge k + 1 \ge 4$ and $C$ is an induced cycle, 
we have 
$E_{G}(u^{+m}, B) \neq \emptyset$, 
and so let $b$ be a vertex of $B$ such that $u^{+m}b \in E(G)$. 
Since 
$m \ge 2$ by Claim~\ref{claim:|C| >= 5}, 
it follows from the definition of $u$ that $b \neq x$. 
These imply that 
(i) or (ii) of Claim~\ref{claim:neighbor of u^{+m}} holds for 
the graph $B$ and the vertices $b, x, u$, 
but 
$u^{+m}b \in E_{G}(\{u^{+m}, u^{-m}\}, V(G-C) \setminus V(B-b))$, 
a contradiction. 
\qed

By Claim~\ref{claim:G - V(C) is not 2-connected}, 
$G-V(C)$ has two end blocks $B_{1}$ and $B_{2}$ 
with cut vertices $b_{1}$ and $b_{2}$, respectively.

\begin{claim}
\label{claim:x_{B}} 
There is a vertex $x_{i}$ 
of $B_{i}-b_{i}$ such that $e_{G}(x_{i}, C) = 2$
for $i \in \{1,2\}$.
\end{claim}

\proof 
Let $B=B_{i}$ and $b_{i}=b$ and 
suppose that every vertex of $B-b$ is adjacent to at most one vertex of $C$. 
Since 
$G$ is $2$-connected, 
there is a vertex $u$ of $C$ such that $E_{G}(u, B-b) \neq \emptyset$, 
i.e., 
some vertex $x$ of $B - b$ is adjacent to $u$. 
Hence 
(i) of Claim~\ref{claim:neighbor of u^{+m}} holds for the graph $B$ and the vertices $b, x, u$, 
and so 
we have 
$E_{G}(u^{+m}, V(G-C) \setminus V(B-b)) = \emptyset$, 
that is, 
$N_{G}(u^{+m}) \subseteq V(C) \cup V(B-b)$. 
Since 
$\delta(G) \ge 4$ and $C$ is an induced cycle, 
this in particular implies that 
$E_{G}(u^{+m}, B-b) \neq \emptyset$.
Then 
again by applying (i) of Claim~\ref{claim:neighbor of u^{+m}} 
for the vertex $u^{+m}$ and a neighbor of $u^{+m}$ in $B-b$ instead of $u$ and $x$, respectively, 
we have 
$E_{G}(u^{+2m}, V(G-C) \setminus V(B-b)) = \emptyset$ 
and thus 
$N_{G}(u^{+2m}) \subseteq V(C) \cup V(B-b)$. 
Repeating this argument 
we get $N_{G}(C) \subseteq V(B-b)$ 
since $|V(C)|$ is odd. 
This implies that $b$ is a cut vertex of $G$, 
which contradicts the $2$-connectivity of $G$. 
Thus 
$e_{G}(x, C) \ge 2$ for some $x \in V(B-b)$. 
Since $C$ satisfies (2) in Lemma~\ref{LM2018 Lemma 5.1}, 
we have $e_{G}(x, C) = 2$. 
\qed

By Claim~\ref{claim:x_{B}} 
and since $C$ satisfies {\rm (2)} in Lemma~\ref{LM2018 Lemma 5.1}, 
each $B_{i} - b_{i}$ contains a vertex $x_{i}$ such that 
$x_{i}u_{i}^{+}, x_{i}u_{i}^{-} \in E(G)$ 
for some $u_{i} \in V(C)$. 

Assume for the moment that 
$u_{1} = u_{2}$. 
Then Claim~\ref{claim:neighbor of u^{+m}} yields that 
$E_{G}(u_{1}^{+m}, V(G-C) \setminus V(B_{1}-b_{1})) = \emptyset$. 
On the other hand, 
since $u_{1}^{+m} = u_{2}^{+m}$, 
it also follows from Claim~\ref{claim:neighbor of u^{+m}} that 
$E_{G}(u_{1}^{+m}, V(G-C) \setminus V(B_{2}-b_{2})) = \emptyset$. 
Therefore we get 
$E_{G}(u_{1}^{+m}, V(G-C)) = \emptyset$. 
But,
since $C$ is an induced cycle, 
this implies that 
$\deg_{G}(u_{1}^{+m}) = \deg_{C}(u_{1}^{+m}) = 2 < 4$, a contradiction. 
Thus 
$u_{1} \neq u_{2}$.

Assume next that 
$u_{1}^{+} = u_{2}$. 
Then $u_{1}^{-m} = u_{2}^{+m}$. 
By using Claim~\ref{claim:neighbor of u^{+m}} 
and arguing as in the above, 
we have $E_{G}(u_{1}^{-m}, V(G-C)) = \emptyset$, which contradicts that $\deg_{G}(u_{1}^{-m}) \ge 4$. 
Thus $u_{1}^{+} \neq u_{2}$. 
Similarly, we have 
$u_{1}^{-} \neq u_{2}$. 
Consequently, 
we get 
\begin{align*}
u_{2} \in V(u_{1}^{+2} \ora{C} u_{1}^{-2}). 
\end{align*}
Since 
$|V(C)| \ge 5$, 
without loss of generality, 
we may assume that 
\begin{align*}
u_{1}^{-} \neq u_{2}^{+}. 
\end{align*}

We will show that
there exist
$2l - 3+\varphi$ $(x_{1}, x_{2})$-paths $\ora{P_{1}}, \dots, \ora{P_{2l- 3+\varphi}}$ 
in $G - V(C)$ 
satisfying 
the length condition. 
In order to show it,
we divide the proof into two cases.

\begin{enumerate}[{\textup{{\bf Case~\arabic{enumi}.}}}]
\setcounter{enumi}{0}
\item 
$\varphi(k) = 1$. 
\end{enumerate}

Since 
$\delta(B_{i}, x_{i}, b_{i}) \ge 2(l-1) + \varphi = 2l-1$ for $i \in \{1, 2\}$, 
it follows from Fact \ref{B is 2-connected} and Theorem~\ref{improvement of LM2018 Lemma 3.1 No.2} 
that 
$B_{1}$ contains $l$ $(x_{1}, b_{1})$-paths $\ora{Q_{1}}, \dots, \ora{Q_{l}}$ satisfying 
the length condition or the semi-length condition,
and 
$B_{2}$ contains $l$ $(b_{2}, x_{2})$-paths $\ora{R_{1}}, \dots, \ora{R_{l}}$
satisfying 
the length condition or the semi-length condition.

Now, 
suppose that either $Q_{1}, \dots, Q_{l}$ or $R_{1}, \dots, R_{l}$ satisfy the length condition.
By the symmetry, 
we may assume that 
$Q_{1}, \dots, Q_{l}$ satisfy the length condition. 
By applying Theorem~\ref{improvement of LM2018 Lemma 3.1} to $(B_{2}, b_{2}, x_{2})$,
we can take other $l-1$ $(b_{2}, x_{2})$-paths $\ora{R_{1}'}, \dots, \ora{R_{l-1}'}$ 
satisfying the length condition, since $\delta(B_{2}, b_{2}, x_{2}) \ge 2(l-1)$. 
Concatenating $Q_{1}, \dots, Q_{l}$ and $R_{1}', \dots, R_{l-1}'$ 
with a $(b_{1}, b_{2})$-path in $G-V(C)$, 
we can obtain 
$2l - 2 \ (=2l - 3+\varphi)$ $(x_{1}, x_{2})$-paths $\ora{P_{1}}, \dots, \ora{P_{2l-2}}$ 
in $G - V(C)$ 
satisfying 
the length condition. 
Therefore, 
we may assume that 
neither $Q_{1}, \dots, Q_{l}$ nor $R_{1}, \dots, R_{l}$ satisfy the length condition, 
and then both of $Q_{1}, \dots, Q_{l}$ and $R_{1}, \dots, R_{l}$
satisfy the semi-length condition.

Let $q$ and $r$ be the switches of 
$Q_{1}, \dots, Q_{l}$ and $R_{1}, \dots, R_{l}$, 
respectively, 
and let $P$ be a $(b_{1}, b_{2})$-path in $G-V(C)$. 
Then by considering the paths
\begin{center}
$x_{1}\ora{Q_{1}}b_{1} \ora{P} b_{2} \ora{R_{i}} x_{2} \ (1 \le i \le r), \ \ \ 
x_{1}\ora{Q_{i}}b_{1} \ora{P} b_{2} \ora{R_{r}} x_{2} \ (2 \le i \le q),$ \\[1mm]
$x_{1}\ora{Q_{q+1}}b_{1} \ora{P} b_{2} \ora{R_{r+1}} x_{2}, \ \ \ 
x_{1}\ora{Q_{q+1}}b_{1} \ora{P} b_{2} \ora{R_{i}} x_{2} \ (r+2 \le i \le l), \ \ \ 
x_{1}\ora{Q_{i}}b_{1} \ora{P} b_{2} \ora{R_{l}} x_{2} \ (q+2 \le i \le l),$
\end{center}
we can obtain $2l - 2 \ (=2l - 3+\varphi)$ $(x_{1}, x_{2})$-paths $\ora{P_{1}}, \dots, \ora{P_{2l-2}}$ 
in $G - V(C)$ 
satisfying 
the length condition.

\begin{enumerate}[{\textup{{\bf Case~\arabic{enumi}.}}}]
\setcounter{enumi}{1}
\item 
$\varphi(k) = 0$. 
\end{enumerate}

Since 
$\delta(B_{i}, x_{i}, b_{i}) \ge 2(l-1)$ for $i \in \{1, 2\}$, 
it follows from Fact~\ref{B is 2-connected} and
Theorem~\ref{improvement of LM2018 Lemma 3.1} that 
$B_{1}$ contains $l-1$ $(x_{1}, b_{1})$-paths 
satisfying 
the length condition 
and 
$B_{2}$ contains $l-1$ $(b_{2}, x_{2})$-paths 
satisfying 
the length condition. 
Concatenating them 
with a $(b_{1}, b_{2})$-path in $G-V(C)$, 
we can obtain 
$2l - 3 \ (=2l - 3+\varphi)$ 
$(x_{1}, x_{2})$-paths $\ora{P_{1}}, \dots, \ora{P_{2l-3}}$ 
in $G - V(C)$ 
satisfying 
the length condition. 

\bigskip
Thus, in both cases,
Table~\ref{k cycles satisfying the length condition} 
gives $2l - 1+ \varphi = k$ cycles 
satisfying the length condition (see also Figure~\ref{lengthconditioncycles}). 

This completes the proof of Theorem~\ref{improvement of LM2018 Theorem 5.2}.  
\qed

\begin{figure}[H]
\begin{minipage}{0.5\hsize}
\begin{center}
\footnotesize
\renewcommand{\arraystretch}{1.5}
\makeatletter
\def\@captype{table}
\makeatother
\begin{tabularx}{60mm}{C}
\hline
$P_{1}, \dots, P_{2l-3+\varphi}$ satisfy the length condition\\
\hline\hline
$x_{1} \, \ora{P_{1}} \, x_{2} \, u_{2}^{-} \, \ola{C} \, u_{1}^{+} \, x_{1}$ \\
$x_{1} \, \ora{P_{2}} \, x_{2} \, u_{2}^{-} \, \ola{C} \, u_{1}^{+} \, x_{1}$ \\
$\vdots$ \\
$x_{1} \, \ora{P_{2l-3+\varphi}} \, x_{2} \, u_{2}^{-} \, \ola{C} \, u_{1}^{+} \, x_{1}$ \\
\hline 
$2l - 3+\varphi$ cycles\\
\hline\\[-6mm]
\hline
$x_{1}  \, \ora{P_{2l-3+\varphi}} \, x_{2} \, u_{2}^{-} \, \ola{C} \, u_{1}^{-} \, x_{1}$ \\ 
$x_{1} \, \ora{P_{2l-3+\varphi}} \, x_{2} \, u_{2}^{+} \, \ola{C} \, u_{1}^{-} \, x_{1}$ \\ 
\hline 
$2$ cycles \\
\hline
\end{tabularx}
\caption{}
\label{k cycles satisfying the length condition}
\end{center}
\end{minipage}
\begin{minipage}{0.5\hsize}
\begin{center}
\input{lengthconditioncycles.tex}
\caption{The path $P_{i}$}
\label{lengthconditioncycles}
\end{center}
\end{minipage}
\end{figure}



\section{Concluding remarks}
\label{sec:concluding remarks}

In this paper, 
we have shown that 
the Thomassen's conjecture on 
the existence of cycles of any length modulo a given integer $k$ 
(Conjecture~\ref{T1983 all length k}) is true 
by giving degree conditions for 
the existence of a specified number of cycles 
whose lengths differ by one or two 
(Theorem~\ref{Conjecture 6.2 of LM2018 for 2-connected case}).

The complete graph of order $k + 1$, in a sense, shows 
the sharpness of the lower bound on the minimum degree condition in Theorem~\ref{Conjecture 6.2 of LM2018 for 2-connected case}. 
On the other hand, 
we believe that 
the assumption of $2$-connectivity in Theorem~\ref{Conjecture 6.2 of LM2018 for 2-connected case} is not necessary. 
In fact, 
Liu and Ma conjectured that 
Theorem~\ref{Conjecture 6.2 of LM2018 for 2-connected case} also holds even 
if we drop the connectivity condition 
(see \cite[Conjecture~6.2]{LM2018}). 
To approach the conjecture, 
the following improvements of 
Theorems~\ref{improvement of LM2018 Lemma 3.1} and \ref{improvement of LM2018 Lemma 3.1 No.2} will be helpful.

\begin{problem}
\label{improvement of LM2018 Lemma 3.1 ver2}
Let $k$ be a positive integer, 
and let 
$(G, x, y)$ be a $2$-connected rooted graph such that $|V(G)| \ge 4$, 
and $z \in V(G)$ $($possibly $z=x$ or $z=y)$. 
Suppose that 
$\deg_{G}(v) \ge 2k$ for any $v \in V(G) \setminus \{x, y, z\}$. 
Then 
$G$ contains $k$ paths from $x$ to $y$ 
satisfying the length condition.  
\end{problem}

\begin{problem}
\label{improvement of LM2018 Lemma 3.1 No.2 ver2}
Let $k$ be a positive integer, 
and let 
$(G, x, y)$ be a $2$-connected rooted graph such that $|V(G)| \ge 4$, 
and $z \in V(G)$ $($possibly $z=x$ or $z=y)$. 
Suppose that 
$\deg_{G}(v) \ge 2k-1$ for any $v \in V(G) \setminus \{x, y, z\}$. 
Then 
$G$ contains $k$ paths from $x$ to $y$ 
satisfying the length condition or the semi-length condition.  
\end{problem}

In fact, 
if Problems~\ref{improvement of LM2018 Lemma 3.1 ver2} and \ref{improvement of LM2018 Lemma 3.1 No.2 ver2} are true, 
then 
by applying them to an end block $B$ with cut vertex $z$ in a given graph of minimum degree at least $k + 1$, 
and 
by arguing as in the proofs of Theorems~\ref{improvement of LM2018 Lemma 4.1} and \ref{improvement of LM2018 Theorem 5.2}, 
we can show that 
$B$ contains $k$ cycles satisfying the length condition when $B$ is $2$-connected but not $3$-connected; 
$B$ contains $k$ cycles, 
which have consecutive lengths or satisfy the length condition when $B$ is $3$-connected and non-bipartite. 
Therefore, 
Problems~\ref{improvement of LM2018 Lemma 3.1 ver2} and \ref{improvement of LM2018 Lemma 3.1 No.2 ver2} 
leads to the improvement of Theorem~\ref{Conjecture 6.2 of LM2018 for 2-connected case}.

The above problems also concern with another Thomassen's conjecture on 
the existence cycles of any even length modulo a given integer $k$.

\begin{refconjecture}[Thomassen \cite{T1983}]
\label{T1983 all even length k}
For a positive integer $k$, 
every graph of minimum degree at least $k + 1$ 
contains cycles of all even lengths modulo $k$. 
\end{refconjecture}

It is known that this conjecture is true for all even integers $k$ (see \cite[Theorem~1.9]{LM2018}). 
The improvement of Theorem~\ref{Conjecture 6.2 of LM2018 for 2-connected case} 
implies that 
it is also true for all odd integers $k$ 
(note that Theorem~\ref{Conjecture 6.2 of LM2018 for 2-connected case} implies that 
Conjecture~\ref{T1983 all even length k} is true for the case where $k$ is odd and a given graph is $2$-connected).

Problems~\ref{improvement of LM2018 Lemma 3.1 ver2} and \ref{improvement of LM2018 Lemma 3.1 No.2 ver2} 
can be proven 
by arguing as in the proofs of Theorems~\ref{improvement of LM2018 Lemma 3.1} and \ref{improvement of LM2018 Lemma 3.1 No.2} 
and by a tedious case-by-case analysis 
(in fact, we have checked Problem~\ref{improvement of LM2018 Lemma 3.1 ver2} is true in a private discussion, but it is unpublished). 
Therefore, we think that giving a short proof for the above problems 
might be interesting and helpful for future work of this research area.








\begin{thebibliography}{99}

\bibitem{BV1998}
J.A.~Bondy, A.~Vince, 
\textit{Cycles in a graph whose lengths differ by one or two}, 
J. Graph Theory \textbf{27} (1998) 11--15. 


\bibitem{LM2018}
C-H.~Liu, J.~Ma, 
\textit{Cycle lengths and minimum degree of graphs}, 
J. Combin. Theory Ser. B \textbf{128} (2018) 66--95. 


\bibitem{T1983}
C.~Thomassen, 
\textit{Graph decomposition with applications to subdivisions and path systems modulo $k$}, 
J. Graph Theory \textbf{7} (1983) 261--271. 


\end{thebibliography}
\end{document}